\documentclass[12pt,a4paper]{article}
\usepackage[utf8]{inputenc}
\usepackage[english]{babel}
\usepackage{amsmath}
\usepackage{amsfonts} 
\usepackage{amssymb} 
\usepackage{makeidx} 
\usepackage{graphicx} 
\usepackage{longtable} 
\usepackage{footnote} 
\usepackage{array}
\usepackage[left=2cm,right=2cm,top=2cm,bottom=2cm]{geometry} 
\usepackage{parskip} 
\usepackage[linguistics, edges]{forest}
\usetikzlibrary{shapes,positioning,trees}
\tikzset{basic/.style = {draw, rounded corners=2pt, rectangle}} 
\usepackage[angle=90,alignment=c, pos={0.9\paperwidth,0.5\paperheight}]{draftwatermark} 
\SetWatermarkText{Procedure to Prove the Collatz Conjecture: Ramachandra Bhat}
\author{RAMACHANDRA BHAT \\ Bengaluru, INDIA.}
\title{Convergence of Collatz Sequences: \\Procedure to Prove the Collatz Conjecture} 
\date{\today} 
\begin{document}  
\maketitle

\abstract{As Collatz conjecture is still to be proved, a method to arrive at the complete proof is explored here.  Conceptually, the process relies on the pre-proven sequence data and the method follows the confirmation of the convergence of the Collatz sequence for all the natural numbers in a sequential forward manner (ascending order of n).  Hence, all the numbers less than n have been proved to converge before beginning to prove the convergence of the sequence for the number n.  

Initially, the the nature of problem was explored and explained the reason for some of the remarkable properties.  Then, the pattern analysis of a first few sequences data led to the concept of pre-proven sequence.  Then, the origin of quartets and subsequently, the octet patterns are described.  The ``Even-Odd" combinations and the concept of (8*v + c) gave a clear picture for the overall process of evaluation for the convergence.  Then onwards, individual set of numbers are analyzed and the convergence confirmations are arrived at.  The process of computing involves the concept of combing.  Finally, the compilation of all the sets together would give the completion of the proof.  Here, the procedure for the computation is presented with example and the current status of computation ($>$99\% completed) while the laborious computation is in progress (for the remaining $<$ 1\%). It is expected that the length of the PSO number list is finite and all the numbers would converge. }

\section*{INTRODUCTION}  

\textbf{Definition:} \textit{The Collatz sequence is defined as the list of numbers that are obtained by the repeat application of the Collatz function on the result of the previous application of the function on a given positive integer, till the number reaches a value of 1(one) wherein the Collatz function is defined as:}
\begin{equation} 
 Col(n) = f(n) = 
\begin{cases}
(3n+1) & \text{if n is odd}\\ (\frac{n}{2}) & \text{if n is even}
\end{cases}
\end{equation} 

\begin{center} 
%\vspace{0.4cm}
\fbox{ \textbf{\emph{Collatz Conjecture$:$}}
\textit{Every natural number (n) generates a Collatz sequence}.} 
\end{center}

\emph{The complete proof of the Collatz conjecture must confirm the convergence of the Collatz sequences for all the natural numbers n = \{1, 2, 3, ..., all positive integers\}}. 

\emph{Total Sequence Order (TSO):} It is the total number of terms present in the trajectory (list of numbers) of the Collatz sequence for the  given natural number (n), i.e., $Col^m(n)$ = 1, where TSO = m, the number of times the Collatz function is applied on n till the sequence reaches the value of 1 for the first instance, i.e., $Col^m(n)$ = 1, where m is the TSO for the given natural number, n.  

\emph{Alternate names$:$} The Collatz conjecture is also known by various alternate names such as (i) Syracuse problem, (ii) 3n+1 (or 3x+1) conjecture, (iii) Hasse’s problem, (iv) Kakutani problem, (v) Thwaites conjecture and (vi) Ulam conjecture.

\section{Background:} 

In the history of mathematics, there are many problems defined through very simple statements and yet, couldn't get to the final solutions due to various reasons.  In literature we could find many such problems listed as ``unsolved problems in Mathematics".  In this era of information, we could get those information through internet searches.   Also available are the information about the enormous amount of effort that are put in by various mathematicians as well as mathematics-lovers to find solutions for these unsolved problems.  The quantum of information available indicates the level of importance and interest as well as the level of  difficulty in solving such problems. 

One of those problems, attempted here to find the solution, is the conjecture that was proposed by Prof. Lothar Collatz in 1937.  Several individuals and organizations have put in a lot of efforts to find the solution to the problem and prove the conjecture, though every explorer understands that it is true but is very difficult to prove.  For any given number, anyone could compute the complete sequence and demonstrate that the sequence converges to a repeat cycle of 4-2-1-4.  However, as a general case, no one has yet proved the convergence for all natural numbers \{N\}.   

It may be recalled from the presentation note by Prof. Terence Tao that Collatz conjecture is still to be proved\cite{Tao}. 
 Hence, an effort was made to understand the problem with a new perspective.  It has been reported in the note that, by computation, the conjecture has been evaluated individually to be true for all the natural numbers from 1 up to $10^{20}$ (confirmed with distributed computing project, 2017).  Out of curiosity to understand the real issue in the proof, an approach to find a generalized procedure for the convergence of a given number is explored here.  

\subsection{Computing a Collatz Sequence(Example)} 

Let us first compute a sequence to understand the process involved and the way sequence appears.  Let us consider n = 7 as a given \textbf{positive odd integer} to demonstrate the same. 

$\therefore$ \emph{Given:} n = 7, an odd positive integer. 

$\therefore$ The first term of the sequence is $Col^1(n)$ = (3n + 1), i.e., $Col^1(7)$ = (3 * 7 + 1) = 22, where $Col(n)$ is the notation used to refer to the term in the sequence of n and 1 refers to first term.  

Since the first term is an even integer, the second term is $Col^2(n)$ = \{Col[$Col^1(n)$]\} = ($\frac{Col^1(n)}{2}$), i.e., $Col^2(7)$ = ($\frac{22}{2}$) = 11. 

Similarly, when we compute the complete sequence for n = 7, the complete trajectory of numbers obtained is as follows: 

\begin{center}
\emph{The trajectory of \textbf{{the Collatz Sequence}} for the given number n =} 7 is
\fbox{Col(n) = Col(7) = 22, 11, 34, 17, 52, 26, 13, 40, 20, 10, 5, 16, 8, 4, 2, 1. } 
\end{center}

Here, the sequence consists of 16 terms and the $16^{th}$ term is 1(one).   It is represented as $Col^{16}(7)$ = 1 and this serial number of the term (=16) is called the total order of the sequence (TSO).    

If we continue to compute the subsequent terms for the sequence, the numbers 4, 2, 1 keeps repeating.  Hence, it is said to reach the repeat cycle of 4-2-1.  However, it is required to show that the sequence reaches the value of 1 (till it reaches for the first time), to prove the convergence of the sequence for the given number.  Also, \textbf{TSO} refers to the sequence reaching the value of 1 for the first time only.  

As we could see here, when an individual natural number is computed for its Collatz sequence, it converges to the 4-2-1 repeat cycle, i.e., when the value of n is known, the sequence computation leads to the repeat cycle.  This is because it is easier and clear to know the nature of the number as even or odd based on its value and then calculate the next term in the series, i.e., to apply the right function ($\frac{n}{2}$) or (3n+1).  For the same reason, there are various \emph{Collatz Sequence Calculator}s available on the internet, as \{n MOD 2\} function could be used to evaluate the nature of the number as even or odd, and compute the complete sequence.    

\subsection{Difficulty in Providing Proof} 

As we could see from the previous section, we could easily compute the sequence for any given even or odd positive integer of known value.  However, when the conjecture needs to be proved, it is necessary to confirm the convergence of all the natural numbers, without any exception, unless otherwise, we could prove the presence of even a single counter-example.  To achieve this, we have to use variables, instead of integers of known value and prove the convergence.  

\subsubsection{Undecidable} 

If the sequence calculation is done with a variable, it would be difficult and may not be possible to judge exactly the value of the number is even or odd, when either possibility exists.  

For example, let us consider a number as \textbf{n = x, as an even positive integer}.  Then, the first term of the sequence is ($\frac{x}{2}$).  We cannot proceed further with the computation unless we know whether this resultant term is even or odd.   

Similarly, let us consider a number as \textbf{n = x, as an odd positive integer}.  Then, the first term of the sequence is (3x+1). Since x is an odd number, this term is an even number.  Hence, the second term of the sequence is ($\frac{3x+1}{2}$).  Now, we cannot proceed further to compute the third term unless we know whether the second term is even or odd.  

That is to say, always any odd number gives the even number as the next number in the sequence whereas the even number when divided by 2 can be either even or odd, depending on the value of the variable in the expression, i.e., if n is odd, then (3n+1) is even, whereas if n is even, then ($\frac{n}{2}$) = ($\frac{even}{2}$) could be either even or odd and hence, undecidable.  This limits the further computation of the sequence.  

While we face the problem of undecidability at the second or the third term of the sequence if we use a variable, same is not true if we compute the sequence for a number of known numberical value.  Hence, we could only compute sequence for any given number of known numerical value.  This gives us the confidence that the conjecture is true and should be provable.  However, the primary difficulty in providing the proof for the conjecture is the undecidability of the natue of the term in the sequence, as even or odd. 

\subsubsection{Unstoppable} 

In addition to this difficulty of undecidability about the value of the term as even or odd, another difficulty we face is to decide when the expression reaches to a value 1(one) where the value of the variable (n) is not known.  While computing the sequence, we always get an expression of the form Col(n) = $\left(\frac{c_1*n + c_2}{2^a}\right)$, where \{$c_1$, $c_2$ and a \} are non-zero positive integers and cannot decide when the value of the term in the sequence becomes 1 (one) to stop the computation\footnote{The computation involves multiplication (by 3), division (by 2) and addition (of 1), but no subtraction.  Hence, $c_1$ and $c_2$ are always non-zero whereas `a' could be zero only for the first few terms in the sequence.}.  One could argue that when the last four terms of the sequence keeps repeating themselves, we could only guess that it has reached the 4-2-1-4 repeat cycle.  However, we won't be able to compute even five terms of the sequence for any number and hence, cannot reach such a stage when the sequences normally contain many terms.   

\subsubsection{Seemingly Impossible} 

One of the methods to attempt to prove the convergence of any series or sequences is the method of Mathematical Induction.  For this, we should be able to calculate a few consecutive number examples to demonstrate the definition of the problem and then to demonstrate the convergence of the series of two consecutive numbers using variables.  In the present scenario of Collatz sequence, though the first part could be done for any integer directly, the second part of it requiring general term of the sequence is difficult (or seemingly impossible) to progress further.  This is because, to write the general term of a sequence we should be able to understand the pattern the terms in the sequence follow.  In the present case, the sequence does not appear to follow any specific pattern and also, the trajectory of the sequence of one number does not give any indication of the trajectory of sequence of any other number directly.  This is the main concern in writing a general expression for the $n^{th}$ term in the sequence and hence, the difficulty in utilizing the concept of \emph{Mathematical Induction} for providing the proof. 

At this stage, for any given number, we have to compute the total sequence for that number independently to understand and confirm the convergence and appears to be impossible to find a way of proving the conjecture using variables.   Hence, the dictum of Paul Erd\H os on the problem \textbf{\lq\lq Mathematics is not yet ready for such problems\rq\rq} has been valid till the publication of this complete work. 

Because of all these difficulties of (i) undecidability, (ii) unstoppability as well as (iii) the seemingly impossible nature of the problem, people have explored to compute the sequence using computers, just to find if there could even be any counter-example to prove or disprove the conjecture by direct computation.  \emph{Direct computation of Collatz sequence for any individual integers cannot provide the proof for the convergence for all the integers, in general.}      

\section{Present Study} 

With the motivational information that (i) there appears to be no known  counter-example as on today, and (ii) all the integers confirm the direct convergence when the sequence is computed individually, the present study was undertaken to explore the solution.  In the present study to understand and overcome these here-mentioned difficulties and to  prove the convergence of the Collatz sequence of any given natural number (n) to the repeat cycle of 4-2-1, the following strategic approach is considered.  

\begin{enumerate}

\item{Compute and compile the sequence for a first few natural numbers, in the sequential ascending order of numbers, to get a pattern of variation of the sequence as we progress by observing the total order of the Collatz sequence (TSO) for these integers.} 

\item{Explore the pattern among the sequences for the numbers such as even or odd integers, multiples of any numbers, exact matching patterns, parallel patterns, etc..  Observe the special behaviour / properties of sequences and explore explanation.  This could bring in clarity for the complete proof.} 

\item{Compute with variables to cover all the integers which are normally difficult to comprehend or compute manually or even using machines (computers).} 

\item{Based on these computations and pattern analysis, define a strategy to progress further and complete the proof.}   

\end{enumerate}
 
\section{Amazing Properties of Collatz Sequences}

Various previous explorers have observed many remarkable properties of the Collatz sequences.  In the present study, some of them have been analyzed and reasoned here.

\begin{enumerate}

\item{An odd number is followed by an even number.  By definition, when n is an odd number, then the next term is (3n+1) = (odd * odd + odd) = (odd + odd) = even.  Hence, Every odd number has to be followed by an even number.} 

\item{Two numbers share the same sequence.} 

\begin{itemize}

For example, 5 and 32 share the same sequence.  Similarly, there are many (infinite) such pairs.  For every odd number, there is an even number which shares the exactly same sequence of terms whereas the vice versa is not true.

Let $n_1$ be an odd integer.  Its first term in the Collatz sequence is (3 * $n_1$) + 1. \\
Let $n_2$ be an even integer. Its first term in the Collatz sequence is $(\frac{n_2}{2})$. \\
If $\frac{n_2}{2} = (3*n_1+1),$ then, the trajectories of the Collatz sequences of both the numbers are exactly same and $n_2 = (6*n_1 + 2)$.   

Hence, the pairs $(n_1, n_2) = (n_1, 6n_1$+2) share the same exact sequence.  For example, (1, 8), (3, 20), (5, 32), (7, 44), (9, 56), etc. Here, it should be noted that $n_1$ is always odd and $n_2$ is always even and $n_1$ is always less than $n_2$.

\end{itemize}

\item{Two consecutive numbers of the form $n_a$=8x+4 and $n_b$=8x+5, have the same total sequence order (TSO) and have exactly same sequence from third term onwards, for all $x>0$}

\begin{itemize}
Let $n_a$ = 8x + 4 and $n_b$ = 8x + 5.  Then, 
\begin{center} 
\begin{tabular}{lcl}

$Col^1(n_a)$ = 4x + 2 & & $Col^1(n_b)$ = 24x + 16 \\
$Col^2(n_a)$ = 2x + 1 & and & $Col^2(n_b)$ = 12x + 8 \\ 
$Col^3(n_a)$ = 6x + 4 & & $Col^3(n_b)$ = 6x + 4 \\

\end{tabular} 
\end{center}
for all integer values of $x > 0$. 
From the third term onwards, the trajectories of both the sequences are exactly same, for the same value of x, for any $x>0$.  

For x=0, it is a special case wherein n=4 already converges to 1 at the second term and the repeat cycle of 4-2-1 starts if computation is continued whereas for n=5, the sequence converges to 1, as the fifth term.  
\end{itemize} 

\item{Two numbers of the form $n_a$=8x and $n_b$=8x+2, have the same total sequence order and have exactly same sequence from fourth term onwards, for all odd values of $x > 1$.}

\begin{itemize}
Let $n_a$ = 8x and $n_b$ = 8x + 2.  Then, 
\begin{center} 
\begin{tabular}{lcl}

$Col^1(n_a)$ = 4x  & & $Col^1(n_b)$ = 4x + 1 \\
$Col^2(n_a)$ = 2x  & and & $Col^2(n_b)$ = 12x + 4 \\ 
$Col^3(n_a)$ = x  & & $Col^3(n_b)$ = 6x + 2 \\
$Col^4(n_a)$ = 3x + 1 (if x:odd) & & $Col^4(n_b)$ = 3x + 1 \\
\end{tabular} 
\end{center}
for all odd integer values of $x > 1$. 
From the fourth term onwards, the trajectories of both the sequences are exactly same, for the same value of x. 

For x=1, it is a special case wherein n=8 already converges to 1 at the third term and the repeat cycle of 4-2-1 starts if computation is continued whereas for n=10, the sequence converges to 1, as the sixth term.  
\end{itemize} 

\item{Third term matching.} 
\begin{itemize}
Let us consider the odd integer values of n.  \\ Then, the first term of the sequence is $Col^1$(n) = (3n+1), an even value term and hence, the second term is $Col^2$(n)=($\frac{(3n+1)}{2}$).

Hence, if we compute the sequence for $n_1$ = n, then we have already computed the sequence for $n_2$ = [$\frac{(3n+1)}{2}$], as it is same as that of $n_1$ from its third term onwards, i.e., the third term in the sequence of $n_1$ and the first term in the sequence of $n_2$ are same and the trajectory continues there onwards.  Hence, the total sequence order for $n_2$ is two terms less than that for $n_1$, i.e., TSO($n_2$) = TSO($n_1$) - 2.    

This is valid only for all odd values of $n_1$.  This helps in calculating the tuple sets\cite{Peter} of Peter Schorer (2017).  For example, (3,5), (5,8), (7, 11, 17, 26), (15, 23, 35, 53, 80), (31, 47, 71, 107, 161, 242), etc.
\end{itemize} 

\item{The term values in the Collatz sequence never repeat.}
\begin{itemize}
If they repeat, then the sequence will get into a loop of its own and does not converge to 1.  Since all the sequences converge to 1 they do not repeat.  Though this is obvious, the observation is taken as is at the beginning of this work and the final proof of the conjectue only proves the fact that the observation is true.

\end{itemize}
 
\item{Surprising trajectories:}
\begin{itemize} 
{The trajectories of  sequences for a few numbers are unexpectedly too long, e.g., 27, 703, 35655, 84975, 151067, etc. the details of which would be obvious only after the completion of the proof of the conjecture.}

Since every odd number has an even number associated with it having exactly the same trajectory of sequence, those even numbers also follow the same surprising trajectory, obviously. 
\end{itemize}  

\end{enumerate} 

Similarly, there are many more remarkable properties for which the proof will be available in the details of the present study while providing the proof of the conjecture.  

\section{Problem exploration}         

As described in the previous section, it is good to know that there exists a lot of repetition of sequences such as exact match, third term matching, fourth term matching, same length sequences, etc.  However, it is also essential to find some pattern to progress further and to explore the possibility to propose a method to confirm the convergence of Collatz sequences.  In this regard, let us first numerically compute a few sequences.  As it is obvious, let us start exploring from n = 1 onward in a sequential forward (ascending order of n) manner.  
  
\subsection{Direct Numerical Computation} 

The Collatz sequence for the first 200 natural numbers is computed and compiled.  Table 1 gives the full sequence details for numbers 1 to 24. The Table provides the information on the number(n), nature of the number as Even or Odd, full trajectory of its Collatz sequence and the total number of terms in the sequence (TSO). 

\begin{center}
\begin{tabular}{||c|c|p{10.9cm}|c||}
  \multicolumn{4}{c}{\textbf{Table 1: The Collatz sequence for the numbers 1 – 24.}}\\[5pt] 
  \hline 
  \multicolumn{1}{||c|}{n} 
  &\multicolumn{1}{|c|}{Type} 
  &\multicolumn{1}{|c|}{Trajectory of the Collatz Sequence } 
  &\multicolumn{1}{|c||}{TSO} \\  
  \hline 1 & Odd & 4, 2, 1, [4, 2, 1, 4, 2, 1....] & 3 \\
  \hline 2 & Even	& 1.  & 1 \\
  \hline 3 & Odd 	& 10, 5, 16, 8, 4, 2, \{1\}. & 7 \\
  \hline 4 & Even	& 2, \{1\}. & 2 \\
  \hline 5 & Odd	& 16, 8, 4, \{2, 1\}. & 5 \\ 
  \hline 6 & Even   & 3, \{10, 5, 16, 8, 4, 2, 1\}. & 8 \\ 
  \hline 7 & Odd	& 22, 11, 34, 17, 52, 26, 13, 40, 20, 10, 5, \{16, 8, 4, 2, 1\}. & 16 \\
  \hline 8 & Even	& 4, \{2, 1\}. & 3 \\ 
  \hline 9 & Odd & 28, 14, 7, \{further trajectory as in n=7\}.  & 19 \\ 
  \hline 10 & Even & 5, \{further trajectory as in n=5\}. & 6 \\
  \hline 11 & Odd & 34, 17, 52, 26, 13, 40, 20, 10, \{further trajectory as in n=10\}. & 14 \\
  \hline 12 & Even & 6, \{further trajectory as in n=6\}. & 9 \\
  \hline 13 & Odd & 40, 20, 10, \{further trajectory as in n=10\}. & 9 \\ 
  \hline 14 & Even & 7, \{further trajectory as in n=7\}. &	17 \\ 
  \hline 15 & Odd & 46, 23, 70, 35, 106, 53, 160, 80, 40, 20, 10, \{5, 16, 8, 4, 2, 1\}.  & 17 \\
  \hline 16 & Even & 8, \{4, 2, 1\}.  & 4 \\ 
  
 \hline 17 & Odd & 52, 26, 13, \{further trajectory as in n=13\}. & 12 \\ 
 \hline 18 & Even & 9, \{further trajectory as in n=9)\}. & 20 \\ 
 \hline 19 & Odd & 58, 29, 88, 44, 22, 11, \{further trajectory as in n=11\}.  & 20 \\ 
 \hline 20 & Even & 10, \{5, 16, 8, 4, 2, 1\}. & 7 \\ 
 \hline 21 & Odd  & 64, 32, 16, \{8, 4, 2, 1\}. & 7 \\
 \hline 22 & Even & 11, \{further trajectory as in n=11\}. & 15\\ 
 \hline 23 & Odd  & 70, 35, 106, 53, 160, 80, 40, 20, \{10, 5, 16, 8, 4, 2, 1\}. & 15 \\
 \hline 24 & Even & 12, \{further trajectory as in n=12\}. & 10\\ 
  \hline 
\end{tabular} 
\end{center}  

The following are the observations from Table 1.  

\begin{enumerate} 

\item{After the sequence once reaches the term value of 1, further continuation of the computation leads to a repeat cycle of 4-2-1-4 endlessly, as shown in n=1.}

\item{Unique integers: Integer 2 (n = 2) is the only integer to have TSO = 1.  Only for n = 5, TSO(n) = n = 5.  Similarly, for n = 16, TSO(n) = $\sqrt{n}$ = 4. }

\item{If the first term of the sequence is of the form $2^a$, where a is any positive integer, the sequence converges directly without any odd numbers till it reaches 1 (e.g., n = 1, 4, 5, 8, 16 \& 21).} 

\item{It may be observed that most of the sequences have many terms in their trajectory, though the major part of the trajectory is the repetition of the sequences computed previously.  It should be noted here that the computation is sequential forward (ascending order of n) and the convergence of the sequence for earlier numbers are proven before computing the next higher number. }  

\item{For all the even numbers, the first term in the trajectory of their sequences is less than the initial number itself, i.e., $Col^1$(n) = $\frac{n}{2}$ $<$ n for all even values of n.  Hence, the further terms in the trajectory are all computed earlier for the other numbers, as the computation is serially ascending integer values of n.}  

\item{Sixth term in the sequence for n = 3 is less than the initial number, i.e., $Col^6$(3) $<$ 3, whereas for n = 5, it is third term, for n = 7 and 15, it is $11^{th}$ term.  This means, there is \textbf{no linear correlation}.   Again, here also, the further terms in their trajectory are all computed earlier for the  other numbers.  \\ This information is useful to reduce the computation load.  Since the earlier sequence is proven to converge to the term value of 1, though the values of some of the subsequent terms in those sequences may even be higher than the starting number, it is useful to confirm the convergence of the number without further computation.} 

\item{The terms in the sequences, after reaching the previously proven number for convergence, are given within the \{ \} curly braces, just to highlight as examples.  A lot of computation could be saved if we keep observing the number pattern in the earlier computed sequences.} 

\end{enumerate} 

Based on these observations, the subsequent computation was done only for the odd numbers and upto the term value less than the given number($<$n), as the remaining terms in the trajectory of the sequence could be just borrowed from the earlier computations.  This sequential forward evaluation method makes the computation easier and also confirms the non-existence of counter example till the number that is last computed. 

\begin{center} 
\begin{tabular}{||c|p{12cm}|c||}
\multicolumn{3}{c}{\textbf{Table 2: The Collatz sequence for the odd numbers between 24 and 64.}}\\[5pt] 
\hline 
\multicolumn{1}{||c|}{n} 
&\multicolumn{1}{|c|}{Trajectory of Collatz Sequence (upto pre-proven sequence term only)} 
&\multicolumn{1}{|c||}{TSO} \\  
\hline 
25 & 76, 38, 19, \{further trajectory as in n=19\}. & 23 \\ 
\hline 27 & 82, 41, 124, 62, 31, 94, 47, 142, 71, 214, 107, 322, 161, 484, 242, 121, 364, 182, 91, 274, 137, 412, 206, 103, 310, 155, 466, 233, 700, 350, 175, 526, 263, 790, 395, 1186, 593, 1780, 890, 445, 1336, 668, 334, 167, 502, 251, 754, 377, 1132, 566, 283, 850, 425, 1276, 638, 319, 958, 479, 1438, 719, 2158, 1079, 3238, 1619, 4858, 2429, 7288, 3644, 1822, 911, 2734, 1367, 4102, 2051, 6154, 3077, 9232, 4616, 2308, 1154, 577, 1732, 866, 433, 1300, 650, 325, 976, 488, 244, 122, 61, 184, 92, 46, 23, \{70, 35, 106, 53, 160, 80, 40, 20, 10, 5, 16, 8, 4, 2, 1\}. & 111 \\ 
\hline 29 & 88, 44, 22, 11, 34, 17, \{further trajectory as in n=17\}. & 18 \\
\hline 31 & 94, \{trajectory as in n=27, 6th term onwards\}. & 106 \\
\hline 33 & 100, 50, 25. & 26 \\ 
\hline 35 & 106, 53, 160, 80, 40, 20. & 13 \\
\hline 37 & 112, 56, 28. & 21 \\ 
\hline 39 & 118, 59, 178, 89, 268, 134, 67, 202, 101, 304, 152, 76, 38. & 34\\ 
\hline 41 & 124, 62, 31. & 109 \\ 
\hline 43 & 130, 65, 196, 98, 49, 148, 74, 37. & 29 \\ 
\hline 45 & 136, 68, 34. & 16 \\ 
\hline 47 & 142, \{trajectory as in n=27, 8th term onwards\}, 46. & 104 \\
\hline 49 & 148, 74, 37. & 24 \\ 
\hline 51 & 154, 77, 232, 116, 58, 29. & 24 \\
\hline 53 & 160, 80, 40. & 11 \\
\hline 55 & 166, 83, 250, 125, 376, 188, 94, 47. & 112 \\
\hline 57 & 176, 86, 43. & 32 \\
\hline 59 & 178, 89, 268, 134, 67, 202, 101, 304, 152, 76, 38. & 32 \\
\hline 61 & 184, 92, 46. & 19 \\
\hline 63 & 190, 95, 286, 143, 430, 215, 646, 323, 970, 485, 1456, 728, 364, \{hereon trajectory continues as in n=27, 18th term onwards\}, 61. & 107. \\ 
\hline 
\multicolumn{3}{r}{\small Results of \itshape Direct numerical computations.}\\
\end{tabular}  
\end{center} 

Further observations from the data of both the tables (n=1-24 and n=25-63) are listed here.  

\begin{enumerate}

\item{While sequence for most of the numbers reaches a pre-proven sequence term value very early, many sequences have very high TSO values.} 

\item{For all even numbers, the first term itself is the pre-proven sequence term value.  Even for odd numbers, there is a pattern with alternate odd numbers reach the pre-proven sequence term value at the third term itself.} 

\item{This concept of pre-proven sequence term value would be of help in overcoming the \textbf{unstoppable} difficulty in the evaluation of the sequences with variables.} 

\end{enumerate}

\subsection{Concept of Pre-Proven Sequences} 

Let us understand this important concept that could help us in overcomming the difficulty of unstoppability, with an example of a numerically computed sequence.  

It should be noted here that the computation is sequential forward and the convergence of the sequence for earlier numbers are confirmed before computing the next higher number. 

Let us consider the first example, \underline{n=14 (an even number)}: 

Trajectory of the Sequence for \{n = 14\}: \{7,\} \{22, 11, 34, 17, 52, 26, 13, 40, 20, 10, 5,\} \{16, 8, 4,\} \{2,\} \{1\}.

Here, the first term of the sequence is 7.  We stop the computation here because, the sequence for n = 7 is already computed and convergence proven.  If we want the complete trajectory, then we could copy the same from the sequence of n = 7.  On similar lines, we could progress till we reach the term value of 1(one).  Hence, for the complete sequence for n = 14, we need the sequence data of n = 7, after that n = 5, after that n = 4, after that n = 2.  The same stepwise computation data of the sequence by the use of pre-proven sequence information is provided here with each segment separated by the curly braces \{\}.  Since the sequences for n = 7, 5, 4 and 2 are already proven for their convergence before the start of the computation for n = 14, it is easy to prove the convergence of sequence for n = 14.  Also, the amount of computation required is drastically reduced, i.e., 1(one) term instead of 17 terms (TSO(n = 14) is 17).  We could also say that all even numbers are divisible by 2 and hence, all the even numbers should converge once we prove the convergence of all the odd numbers.  Therefore, there should be no necessity to prove the convergence of any of the even numbers.     

Let us consider another example, \underline{n = 59 (an odd number)}.

Trajectory of the Sequence for \{n = 59\}: \{178, 89, 268, 134, 67, 202, 101, 304, 152, 76, 38,\} \{19,\} \{58, 29, 88, 44, 22, 11,\} \{34, 17, 52, 26, 13, 40, 20, 10,\} \{5,\} \{16, 8, 4,\} \{2,\} \{1\}.

For n = 59, the sequence reaches a pre-proven term value of 38 at the term serial number of 11.  We could stop the computation here because, the sequence for n = 38 is already computed and convergence proven.  Similar to the earlier scenario of n = 14, we could now complete trajectory by copying those pre-approved sequence data for n = 38, 19, 11, 10, 5, 4, 2.  These are provided here with each segment separated by the curly braces \{\}.  Hence, the sequence for n = 59 is confirmed to converge merely by computing the sequence only upto the pre-proven sequence term value (11 terms) instead of 32 terms, ($\because$ TSO(n = 59) is 32).  This demonstrates the reduced computation need even for odd integers.  However, it should be noted that the reduction in computation need for odd integers is not uniform as compared to that for the even integers.  Therefore, we have to compute and confirm the convergence of each odd integers to prove the conjecture.  

It is easy to compute the sequence for any given number of known integer value.  But, as described in the unstoppability of the nature of computation for a number with a variable, it is difficult to know the computation termination step.  If we use variable (n) to compute the sequence, we could always compare the value of the computed term in the sequence with that of the variable (n), i.e., \{term value $<$ n ?\}.  With this, we could overcome the difficulty of unstoppability.  Hence, the concept of pre-proven sequence is helpful.  

\begin{center} 
\begin{verse}
Therefore, \emph{Pre-proven sequence concept is defined as the use of the information of already computed and confirmed convergence of the sequence of a number of value lesser than the number for which the sequence is being computed}.
\end{verse} 
\end{center} 

Another term used in this study is the \emph{Pre-proven sequence order (PSO)}.  It is defined as the order (serial number of the term in the trajectory) of the Collatz Sequence at which the term value of the sequence reaches the value less than the initial value (n) for the first instance \{i.e., $Col^{PSO}(n) < n$, at the first instance\}.  On further continuation of the computation, the term values may go even higher than n, but the sequence should eventually reach the value of 1, i.e., TSO.  Since this later part of the sequence is already computed in a different case and confirmed the convergence, we should just use the same data, without re-computing the same sequence trajectory again.   Therefore, if the sequence reaches the PSO term [$Col^{PSO}$(n) $<$ n] during the  computation, the sequence could be considered to converge to 1(one), where PSO is the pre-proven sequence order, whereas TSO is the total sequence order when the sequence actually reaches to 1 at first instance [$Col^{TSO}$(n) = 1].   Hence, it may be noted that, always, PSO $\le$ TSO.

\emph{For example,} the Collatz sequence for the number n=9 is \{28, 14, 7, [22, 11, 34, 17, 52, 26, 13, 40, 20, 10, 5, 16, 8, 4, 2, 1 ]\}.  The total number of terms in the sequence is 19 and hence TSO(n=9)=19.  However, the third term in the sequence, viz., $Col^3(9)$=7, is less than n(=9) and hence, PSO(n=9) = 3 and the proven sequence term value (PST value) is 7. 

Also, to know the TSO value for the number for which the sequence is being computed, once the PSO is reached, we just could use the relation TSO(n) = PSO(n) + TSO$\{Col^{PSO}$(n)\}, where \{$Col^{PSO}$(n)\} is the term value of  $PSO^{th}$ term of the Collatz sequence for n.

As we go along we will understand that \underline{\textbf{this appears to be the only way}} we could attempt to prove the convergence of Collatz sequences and hence, the proof for the Collatz Conjecture.  This way, we overcome the problem of unstoppability.    

To confirm the existence of any pattern and utilize the same information in our study, further computation of the trajectories of the sequences were evaluated only upto the PSO term with the actual numbers upto n=200, as the remaining terms in the trajectory of the sequence could be just borrowed from the earlier computation.  These data are tabulated in Table 3.  Here, only the number(n), PSO and TSO values are tabulated for easy reference.  

\begin{center}
\begin{longtable}[ht]{||c|c|c||c|c|c||c|c|c||c|c|c||} 
  \multicolumn{12}{c} 
     {\textbf{Table 3: The summary of Collatz sequences for the numbers 1-200.}} \\[5pt] 
  \hline 
    \multicolumn{1}{||c|}{n} 
   &\multicolumn{1}{|c|}{PSO} 
   &\multicolumn{1}{|c||}{TSO} 
   &\multicolumn{1}{|c|}{n} 
   &\multicolumn{1}{|c|}{PSO} 
   &\multicolumn{1}{|c||}{TSO} 
   &\multicolumn{1}{|c|}{n} 
   &\multicolumn{1}{|c|}{PSO} 
   &\multicolumn{1}{|c||}{TSO} 
   &\multicolumn{1}{|c|}{n} 
   &\multicolumn{1}{|c|}{PSO} 
   &\multicolumn{1}{|c||}{TSO} \\ 
 \hline 
  \endfirsthead 
     \multicolumn{12}{c} 
       {\textbf{The summary of the sequences for the numbers 1-200.}
         (\textit{Continued})}\\[5pt] 
  \hline 
  \multicolumn{1}{||c|}{n} 
  &\multicolumn{1}{|c|}{PSO} 
  &\multicolumn{1}{|c||}{TSO} 
  &\multicolumn{1}{|c|}{n} 
  &\multicolumn{1}{|c|}{PSO} 
  &\multicolumn{1}{|c||}{TSO}  
  &\multicolumn{1}{|c|}{n} 
  &\multicolumn{1}{|c|}{PSO} 
  &\multicolumn{1}{|c||}{TSO} 
  &\multicolumn{1}{|c|}{n} 
  &\multicolumn{1}{|c|}{PSO} 
  &\multicolumn{1}{|c||}{TSO} \\ 
  \hline 
  \endhead 
  \hline 
  \multicolumn{12}{r}{\small\itshape{continued on the next page}}\\
  \endfoot 
  \hline 
  \multicolumn{12}{r}{\small Direct numerical computations.}
  \endlastfoot  
  \hline - & - & -  &  1 & NA &  3 & 2 & 1 & 1 & 3 & 6 & 7 \\ 
  \hline 4  & 1 & 2  &  5 & 3 &  5 & 6 & 1 & 8 & 7 & 11 & 16 \\ 
  \hline 8  & 1 & 3  &  9 & 3 & 19 & 10 & 1 & 6 & 11 & 8 & 14 \\ 
  \hline 12 & 1 & 9  & 13 & 3 &  9 & 14 & 1 & 17 & 15 & 11 & 17 \\ 
  \hline 16 & 1 & 4  & 17 & 3 & 12 & 18 & 1 & 20 & 19 & 6 & 20 \\
  \hline 20 & 1 & 7  & 21 & 3 &  7 & 22 & 1 & 15 & 23 & 8 & 15 \\
  \hline 24 & 1 & 10 & 25 & 3 & 23 & 26 & 1 & 10 & 27 & 96 & 111 \\
  \hline 28 & 1 & 18 & 29 & 3 & 18 & 30 & 1 & 18 & 31 & 91 & 106 \\
  \hline 32 & 1 & 5  & 33 & 3 & 26 & 34 & 1 & 13 & 35 & 6 & 13 \\
  \hline 36 & 1 & 21 & 37 & 3 & 21 & 38 & 1 & 21 & 39 & 13 & 34 \\
  \hline 40 & 1 & 8  & 41 & 3 & 109& 42 & 1 &  8 & 43 & 8 & 29 \\
  \hline 44 & 1 & 16 & 45 & 3 & 16 & 46 & 1 & 16 & 47 & 88 & 104 \\
  \hline 48 & 1 & 11 & 49 & 3 & 24 & 50 & 1 & 24 & 51 & 6 & 24 \\
  \hline 52 & 1 & 11 & 53 & 3 & 11 & 54 & 1 & 112& 55 & 8 & 112 \\
  \hline 56 & 1 & 19 & 57 & 3 & 32 & 58 & 1 & 19 & 59 & 11 & 32 \\
  \hline 60 & 1 & 19 & 61 & 3 & 19 & 62 & 1 & 107& 63 & 88 & 107 \\
  \hline 64 & 1 &  6 & 65 & 3 & 27 & 66 & 1 & 27 & 67 &  6 & 27 \\
  \hline 68 & 1 & 14 & 69 & 3 & 14 & 70 & 1 & 14 & 71 & 83 & 102 \\
  \hline 72 & 1 & 22 & 73 & 3 & 115& 74 & 1 & 22 & 75 &  8 & 14 \\
  \hline 76 & 1 & 22 & 77 & 3 & 22 & 78 & 1 & 35 & 79 & 13 & 35 \\
  \hline 80 & 1 &  9 & 81 & 3 & 22 & 82 & 1 & 110 & 83 & 6 & 110\\
  \hline 84 & 1 &  9 & 85 & 3 &  9 & 86 & 1 & 30 & 87 &  8 & 30 \\
  \hline 88 & 1 & 17 & 89 & 3 & 30 & 90 & 1 & 17 & 91 & 73 & 92 \\
  \hline 92 & 1 & 17 & 93 & 3 & 17 & 94 & 1 & 105 & 95 & 13 & 105 \\
  \hline 96 & 1 & 12 & 97 & 3 & 118 & 98 & 1 & 25 & 99 & 6 & 25 \\
  \hline 100 & 1 & 25 & 101 & 3 & 25 & 102 & 1 & 25 & 103 & 68 & 87 \\
  \hline 104 & 1 &  12 & 105 & 3 &  38 & 106 & 1 &  12 & 107 &  8 & 100 \\
  \hline 108 & 1 & 113 & 109 & 3 & 113 & 110 & 1 & 113 & 111 & 50 & 69 \\
  \hline 112 & 1 &  20 & 113 & 3 &  12 & 114 & 1 &  33 & 115 &  6 & 33 \\
  \hline 116 & 1 &  20 & 117 & 3 &  20 & 118 & 1 &  33 & 119 &  8 & 33 \\
  \hline 120 & 1 &  20 & 121 & 3 &  95 & 122 & 1 &  20 & 123 & 13 & 46 \\
  \hline 124 & 1 & 108 & 125 & 3 & 108 & 126 & 1 & 108 & 127 & 24 & 46 \\
  \hline 128 & 1 &   7 & 129 & 3 & 121 & 130 & 1 &  28 & 131 &  6 & 28 \\
  \hline 132 & 1 &  28 & 133 & 3 &  28 & 134 & 1 &  28 & 135 & 11 & 28 \\
  \hline 136 & 1 &  15 & 137 & 3 &  80 & 138 & 1 &  15 & 139 &  8 & 41 \\
  \hline 140 & 1 &  15 & 141 & 3 &  15 & 142 & 1 & 103 & 143 & 11 & 103 \\
  \hline 144 & 1 &  23 & 145 & 3 & 116 & 146 & 1 & 116 & 147 &  6 & 116 \\
  \hline 148 & 1 &  23 & 149 & 3 &  23 & 150 & 1 &  15 & 151 &  8 & 15 \\
  \hline 152 & 1 &  23 & 153 & 3 &  36 & 154 & 1 &  23 & 155 & 65 & 85 \\
  \hline 156 & 1 &  36 & 157 & 3 &  36 & 158 & 1 &  36 & 159 & 34 & 54 \\
  \hline 160 & 1 &  10 & 161 & 3 &  98 & 162 & 1 &  23 & 163 &  6 & 23 \\
  \hline 164 & 1 & 111 & 165 & 3 & 111 & 166 & 1 & 111 & 167 & 47 & 67 \\
  \hline 168 & 1 &  10 & 169 & 3 &  49 & 170 & 1 &  10 & 171 &  8 & 124 \\
  \hline 172 & 1 &  31 & 173 & 3 &  31 & 174 & 1 &  31 & 175 & 13 & 80 \\
  \hline 176 & 1 &  18 & 177 & 3 &  31 & 178 & 1 &  31 & 179 &  6 & 31 \\
  \hline 180 & 1 &  18 & 181 & 3 &  18 & 182 & 1 &  93 & 183 &  8 & 93 \\ 
  \hline 184 & 1 &  18 & 185 & 3 &  44 & 186 & 1 &  18 & 187 & 11 & 44 \\
  \hline 188 & 1 & 106 & 189 & 3 & 106 & 190 & 1 & 106 & 191 & 21 & 44 \\
  \hline 192 & 1 &  13 & 193 & 3 & 119 & 194 & 1 & 119 & 195 &  6 & 119 \\
  \hline 196 & 1 &  26 & 197 & 3 &  26 & 198 & 1 &  26 & 199 & 13 & 119 \\
  \hline 200 & 1 &  26 & 201 & 3 &  18 & 202 & 1 &  26 & 203 &  8 &  26 \\ 
  \hline 
\end{longtable} 
\end{center} 

Let us summarize all the observations here from all the data generated so far (trajectories of Collatz sequences for n = 1-200): 

\begin{enumerate}

\item{The full Collatz sequences for the first three integers, viz., n = 1, 2 and 3 are computed directly as the initial data.  n = 2 is the only number with a single term in its sequence.  n = 3 has 7 terms in its trajectory with 16 as the highest individual term value whereas for n = 27, TSO(27) = 111, PSO(27) = 96 and the highest individual term value is 9232. } 

\item{For all even numbers, the first term itself is the pre-proven sequence term value, though the total sequence order (TSO) values are varying. }

\item{Even for odd numbers, there is a pattern with alternate odd numbers reach the pre-proven sequence term value at the third term itself.} 

\item{While sequences for many numbers reach a pre-proven sequence term (PSO) value very early, some of the sequences have very high TSO values.} 

\item{Majority of the numbers have the PSO values less than 30.  Even those numbers having high PSO values are having some sort of a pattern.} 

\item{The \textit{concept of pre-proven sequence} term value would seem to be of great help in understanding the pattern of number types and also, to overcome the difficulty of \textbf{unstoppable} nature of computation of sequences during the exploration of the sequences with variables.} 

\end{enumerate}

\subsection{Pattern Analysis} 

The PSO and TSO values for all the numbers from 1 to 200 are tabulated in Table 3.  While tabulating these data, a careful observation led to the classification of the numbers into four category of numbers, based on their PSO values.  While the PSO values give a very pleasant surprise, the combination of PSO and TSO values also give an insightful picture.  The overall observations are as follows: 

\begin{enumerate}

\item{It should be noted here that the computation is sequential and the convergence of the sequence for earlier numbers are proven before computing the next higher number.}

\item{There is some sort of a pattern in the PSO and TSO values and only a few ($<$ 3\% numbers) have high PSO values ($>$ 32).  Even for numbers exhibiting high TSO values, many of them exhibit low PSO values only.}  

\item{There are many sets of two consecutive numbers having the same TSO values.  Also, there are a few pairs of numbers exhibiting exactly same trajectory of sequences.} 

\item{All the even integers exhibit a PSO value of 1(one), as expected obviously by the definition of the Collatz function.} 

\item{Among all the odd positive integers, half of them exhibit a PSO value of 3(three).  Hence, this also reduces the overall computation load. } 

\item{Hence, the computation required to be demonstrated is only for one type of odd integers, effectively reducing the computation load by 75\%.} 

\item{Even among those remaining 25\% integers, we could observe that majority of them converge to a PSO value of less than 15.  Primarily, 6, 8, 11 and 13 are the PSO values observed with higher frequently relative to other values.} 

\item{Even though the numerical computation of the Collatz sequence for the first few numbers gave an important concept of pre-proven sequence, the pattern analysis takes us on a different path, i.e., we cannot progress with any sequentially forward computation with variables.  Instead, we group the natural numbers as different sets and confirm the convergence of the set of numbers to their PSO term.  This leads to rephrasing the conjecture as follows:} 

\end{enumerate} 

\begin{center} 
\begin{verse}
\emph{Every natural number reaches its pre-proven Collatz sequence term} [and hence, generates a Collatz sequence]. That is to say, \textbf{for every natural number n there exists a number k such that $Col^k(n) < n $. }Hence, the revised problem definition is to prove that all the natural numbers reach PSO term.  Then, the concept of sequentially forward computing and the concept of pre-proven sequence together proves the conjecture.    

\end{verse} 
\end{center} 

\section{Procedure for computation} 

Initially, the Collatz sequences for the first three numbers 1, 2 and 3 are computed numerically  and confirmed the convergence. This is the initial data that is required to confirm the convergence of the subsequent numbers that follow during the computation of all the natural numbers, as the computation is sequentially forward\footnote{It should be noted here that the following two phrases, viz., \emph{sequence converges to a PSO value} or \emph{sequence converges to 1},  are used interchangeably to mean the confirmation of convergence as the sequence reaching the PSO term means the convergence of the sequence to a term value of 1(one) subsequently, on the basis of the pre-proven sequence concept.}.

Computing the Collatz sequence for any given number is a straightforward procedure.  At every step, check whether the term is even or odd, apply the appropriate function to derive the next term.  Once the term value reaches 1(one), stop the computation and the sequence is ready.  The total sequence order (TSO) may be found by counting the number of terms present in the trajectory of the sequence.  Similarly, the pre-proven sequence order (PSO) may be found by counting the number of terms present in the trajectory of the sequence till the first occurrence of the term (including that first occurrence  term) of value less than the initial number. 

However, when we \textbf{compute using variables} instead of the integers (natural numbers) to explore the sequence, as we have already described, we come across the problems of (i) \textbf{unstoppable} and (ii) \textbf{undecidable} characteristics of the process.  

\subsection{Stop at the Past} 

Now, we overcome the problem of unstoppability by the concept of pre-proven sequence order.  Once the term computed is of value less than the initial number, the computation is stopped.  Since the computation is sequentially forward, the remaining part of the sequence is already computed earlier for a number of the same value as that of the PSO term.    

\subsection{Decide to be undecided}

While using variables for the computation of the sequence, it is difficult to understand and decide whether the term is of even value or odd value, particularly so when the term is expressed as a fraction (denominator is $2^a$, a$>$0), though the actual numerical value is an integer.  To overcome this problem, it is better and necessary also to consider both the options (even and odd), one after the other, to find the next term and subsequently to the PSO term.  Hence, at every point when we come across the  terms of unknown nature as even or odd, we first consider that the term is an even value term and complete the computation for that option.  Next, we take up the option of the term being the odd value term and compute the sequence for the PSO term.  With this, we could overcome the problem and explore the complete sequence, i.e., till the PSO term.

With these two considerations, we could compute the sequences using variables and find the values of variables for which the sequence converges to its PSO term and hence, would subsequently converge to the term value of 1(one). To confirm that we do not leave any possible options unattended, we use the rolling/folding steps to comb all the options.  During the combing operation, the initially even term-value options are computed till the PSO term and then the odd term-value options are kept adding to the "to-be computed" list and taken subsequently with "first-in first-out" concept for computation.  

\section{Origin of the quartets} 

From the summary of the Collatz sequences for the first 200 natural numbers (cf Table 3), it may be observed that the natural numbers form 4 different sets with respect to the PSO values wherein the two sets exhibit PSO value of 1 and one set exhibit PSO value of 3 and the fourth set has varying values of PSO.   The same observation could be explained as follows: 

Let the natural number be n. 

Then, if n is even, the first term of the sequence is $\left(\cfrac{n}{2}\right)$ which is obviously less than the initial number (n).  Hence, the PSO for all the even numbers is 1(one).  

Therefore, for n to be even, let n = 2y, where y = \{1, 2, 3, ....(all positive integers\footnote{Since Collatz function is defined only for n $>$ 0, zero is not included.})\}.    

If n is odd, the first term of the sequence is (3n+1) which is greater than n.  Since n is odd, this first term is even \{(odd*odd+odd)=even\}.  On further application of the Collatz function, second term is ($\frac{3n+1}{2}$) which is obviously not less than the initial number (n).  Hence, to get to the PSO term, we need to apply the Collatz function again.  

Hence, $Col^1(n)$ = (3n+1). 

$Col^2(n)$ = $\left(\frac{3n+1}{2}\right)$. 
  
$Col^3(n)$ = $\left(\frac{3n+1}{4}\right)$, iff $\left(\frac{3n+1}{2}\right)$ is even.  In this case, the third term is less than the initial number for all n $>$ 1, i.e., PSO = 3.  It should be noted that n = 1 is a special case where the initial number itself is 1(one) and the 4-2-1 cycle repeats on repeated application of the Collatz function.  Also, Since n=1 is the first number in the set of natural numbers, the concept of PSO is not applicable. 
 
Here, n is odd and $\left(\frac{3n+1}{2}\right)$ is even implies, n = \{1, 5, 9, ...(alternate odd numbers)\}.  The same may be written in the form n = (4P+1), wherein P = \{0, 1, 2, 3, ...(all positive integers)\}.  

Similarly, if n is odd and $\left(\frac{3n+1}{2}\right)$ is also odd, then, n = \{3, 7, 11, 15, ... (the remaining alternate odd numbers)\}. These numbers may be written in the form n = (4P+3), wherein P = \{0, 1, 2, 3, ...(all positive integers)\}.  In this case, the computing to get the PSO term continues, as each set of numbers exhibit varying PSO values. 

Since the odd numbers are grouped as (4P+1) and (4P+3), we may re-write the even numbers also as (4P+0) and (4P+2), for uniformity considering options of y=Even=(2P) and y=Odd=(2P+1).  Therefore, the whole number set is written as \{N\} = \{1, 2, 3, (4P+0), (4P+1), (4P+2) and (4P+3)\}, wherein P = \{1, 2, 3, ...(all positive integers)\}.  Hence, the pattern analysis led to the grouping of the numbers as quartets.  This grouping is convenient to evaluate the PSO values for the remaining number sets of the form n = (4P+3), where the PSO values vary.

\begin{center} 
\fbox{ \{N\} = \{1, 2, 3, (4P+0), (4P+1), (4P+2), (4P+3)\}, where P = \{1, 2, 3, .....\} }   
\end{center}  

\section{Quartet Analysis} 

\subsection{EVEN Numbers}  

\subsubsection{The Integers of the form n = 4P.}  

The first term of the sequence for n = 4P is ($\frac{n}{2}$) = 2P.  i.e., $Col^1$(4P) = 2P. \\ This is less than the initial number n(=4P), for all values of P.  $\therefore$  PSO(n=4P) = 1. \\Therefore, the Collatz sequence for all the numbers of the form n = 4P, converges to a PSO value of 1 and subsequently, reaches to a sequence term value of 1(one), for all P. 
 
$\therefore$ PSO(n=4P) = 1 and TSO(n=4P) = 2 + TSO(n=P). 
 
\subsubsection{The Integers of the form n = (4P+2).}  

The first term of the sequence for n = (4P+2) is \{$\frac{n}{2}$\} = (2P+1) $<$ (4P+2).  $\therefore$ PSO\{n = (4P+2)\} = 1.  Hence, the Collatz sequence for all the numbers of the form n = (4P+2), converges to a term value of 1(one) for all P. 

$\therefore$ PSO(n=4P+2) = 1 and TSO(n=4P+2) = 1 + TSO(n=2P+1). \\With this we cover all the even natural numbers.  This is with the assumption that the sequences for all the even and odd numbers up to the number \{n = 4P or n = (4P+2)\} for which the sequence being computed are converging to the 4-2-1 repeat cycle without any exception, i.e., without any counter-example.  

\subsection{ODD Numbers}

\subsubsection{The Integers of the form n = (4P+1).} 

\underline{Sequence}:   

First term:  $Col^1$\{n = (4P+1)\} = {3(4P+1) + 1} = (12P + 4) : (even number).

Second term:  $Col^2$\{n = (4P+1)\} = $Col^1(12P+4)$ = ($\frac{12P+4}{2}$) = (6P+2) : (even number). 

Third term:	 $Col^3$\{n = (4P+1)\} = $Col^1(6P+2)$ = ($\frac{6P+2}{2}$) = (3P+1) : (number is even or odd, depends on whether value of P is even or odd). 

However, since (3P+1) – the third term of the sequence is less than the initial number n \{= (4P+1)\}, it is not essential to know whether that term is odd or even as the computation would be terminated at that stage ($\because$ sequence has already reached the pre-proven sequence term value).  Because the computation is sequentially forward and all the numbers up to n =(4P+1) should have already been proven to converge, all the numbers of the form n=(4P+1) also converges, with a PSO value of 3. 

$\therefore$ PSO(n=4P+1) = 3 and TSO(n=4P+1) = 3 + TSO(n=3P+1). 

With this 75\% of the natural numbers are expected to converge without any counter example, with \emph{a caviet} that all the other remaining (25\%) numbers also converge without any counter-example. 

Hence, the proof of the Collatz Conjecture depends completely on the proof of the existence or absence of the counter-examples among the numbers of the form n=\{4P+3\}, where P is any positive integer (including zero).  

\subsubsection{The Integers of the form n = (4P+3).}  
		
Here, \textbf{the values of n are odd only}, whereas the values of P could be either even or odd, including zero.  This is the case where all the variability and undecidability crops up while computing the sequences (till the PSO term) for any given number.  The difficulty in knowing the value of the term as even or odd depends on the value of P as either even or odd and also the fractions involving P when we divide the even term by 2.  This poses the problem in progressing further in computing the subsequent terms.  Also, looking at the extreme right block of n-PSO-TSO values in Table 3 where the numbers are in the form of n = (4P+3), we could get an idea of the variability and complexity in computing the total sequence for any given number.  However, the good news within those difficult scenarios is that a lot of repeat occurrences for PSO = 6, 8, 11 \& 13 could be observed.  The pattern of repeat occurrences gives us an initial hint that we need to evaluate the sequences separately for even and odd values of P.  With this initial hint, let us consider the computation of sequences for any given value of n = (4P+3) with P as a variable number.  

\subsection{Sequence Computation: n=(4P+3):Rules of the GAME} 

In the previous sections, we have seen that 75\% of all the integers do not contribute to the existence of any counter-examples.   Also, they do not pose any problem in finding the pre-proven sequence order, i.e., there is no difficulty of decision making or stopping the computation.  But, only the odd value integers of the form n = (4P+3) are the only set of integers those exhibit complex computation scenario, as could be understood from Table 3, section 2.  Also, a considerable portion of the set of integers of this form also exhibit PSO values in the range of 6, 8, 11 and 13, though there are other integers exhibiting PSO values more than 100.  Now is the time when we have to face all the difficulties of computation and decision making.  Therefore, let us begin with the computation of sequences for any given value of n = (4P+3) with P as an integer variable, understand the actual difficulties as we stumble upon and resolve them, as appropriate.   

As the intensity of the computation starts becoming complex, it is better to specify the rules of the game, i.e., the conditions for computations that need to be complied with.  

\begin{itemize}

\item{The first two terms of the sequences are of the form (3n+1) and ($\frac{3n+1}{2}$).  

$1^{st}$ term: 
$Col^1(n)$ = $(n_1)$ = {$\underbrace{(3*n + 1)}_{Even\ only}$} = {$\underbrace{(3*n+1)}_{Even}$} +  {$\underbrace{(0*n+0)}_{Even} $}.  

$2^{nd}$ term: 
$Col^2(n)$ = $(n_2)$ = {$\underbrace{\frac{(3*n + 1)}{2}}_{Odd\ or\ Even}$} = {$\underbrace{(1*n+0)}_{Odd}$} + {$\underbrace{\frac{(1*n+1)}{2}}_{Odd\ or\ Even ?}$}.} 

\item{The general form of the terms in the sequences is of the form: 

\vspace{15pt} 

\begin{center} 
\fbox{
$a^{th}$ term: 
$Col^a(n)$ = $(n_a)$ = {$\underbrace{\frac{(b*V + c)}{d}}_{Odd\ or\ Even}$} = {$\underbrace{(e*V+f)}_{Odd\ or\ Even}$} +  {$\underbrace{\frac{(g*V+h)}{d}}_{Odd\ or\ Even} $} }
\end{center} 
\vspace{15pt} 

where V is the variable and n = f(V), the given integer, b=(d*e+g) and c=(d*f+h).  Also, the values of a, b, c, d, e, f, g and h as well as the expressions ($\frac{b*V + c}{d}$) and ($\frac{g*V+h}{d}$) are integers only and no fractions.}  

\item{Invariably, for the first term d=1.  For subsequent terms, d$\ge$ 1.  Once the value of d is greater than 1, the undecidability problem starts.  As a strategy, we do not decide.  Instead of struggling to decide on whether the term is even or odd without the appropriate information for the same, we consider the term value as even and evaluate the sequence and then come back to that same term and evaluate the sequence considering the term value as odd.  Since we consider both the possible options without deciding on either as even or odd, this strategy is referred to as "\textbf{decide to be undecided}" and overcome the problem of undecidability by exploring both the options.  By this method, we do not miss on any number or set of numbers those may even exhibit as counter-examples.}

\item{Since the computation, now onwards, is only for the odd positive integers of the form n = (4P+3), \emph{n is always an odd integer} whereas, in any particular case of computation, P can be of even or odd integer value (including zero).  Also, there are no terms with a value of 0(zero) in all the computations. }

\item{It should be noted here that the following rules of addition, multiplication and division are used in the decision making at various steps during the computations.  One special observation in the process of computation is that there are no steps where subtraction is involved.} 

\begin{itemize} 
\item{Rules of number addition: \\
Even + Even = Even; Odd + Odd = Even;\\  
Even + Odd = Odd; Odd + Even = Odd } \\ 

\item{Rules of number multiplication: \\
Even * Even = Even; Odd * Odd = Odd; \\ 
Even * Odd = Even; Odd * Even = Even}\\

\item{Division possibilities: 

All the computations involved in this study are only with integers and no fractions.  It should be noted while doing the computation that any term in the Collatz sequence should only be a positive integer and no negative integer, or zero, or any fraction.  Only even numbers are divided by a factor of 2, and the quotient should be an integer only. 

i.e., ($\cfrac{Even}{2}$) = Even or Odd, the quotient is integer only and no fraction.}

\end{itemize}

\item{While computing the sequences, we stop the computation by comparing the value of the term with that of the initial number (term value $<$ n ?) to decide whether we have reached the PSO term.  This way, we overcome the problem of unstoppability.  Since we are trying to confirm that the sequence has reached the pre-proven term value that has been proven earlier, this strategy is considered as "\textbf{Stop at the Past}".  

Here onwards, each individual case of P values will be analyzed for the convergence of the Collatz sequences for those individual set of numbers. To evaluate the set of values that satisfy the conditions of convergence, these rules are important.}

\end{itemize}

\subsubsection{Sequence Computation for n = (4P+3).}

Let us begin the computation for \textbf{n = (4P + 3)}, with  P is a variable. 

At every step of computation, we need to assess (i) the value of the term as even or odd and (ii) whether the value of term is less than the value of the initial number [n = (4P+3)].  It is to decide whether to progress the computation or to stop if the term value reaches the pre-proven sequence (PSO) term value.   

$1^{st}$ term = $Col^1$(n) = \{$\underbrace{3*(4P + 3) + 1}$\} = {$\underbrace{(12P+10)}_{Even\ only}$}

$2^{nd}$ term = $Col^2$(n) = \{$\underbrace{\frac{(12P + 10)}{2}}$\} = {$\underbrace{(6P+5)}_{Odd\ only}$}

$3^{rd}$ term = $Col^3$(n) = \{$\underbrace{{3*(6P + 5) + 1}}$\} = {$\underbrace{(18P+16)}_{Even\ only}$}

$4^{th}$ term = $Col^4$(n) = {$\underbrace{\frac{(18P + 16)}{2}}$} = {$\underbrace{(9P+8)}_{Odd\ or\ Even\ ?}$}  

Now, the fourth term $Col^4$(n) could be of either even-value or odd-value depending on whether P is even or odd.  It is even, if P is even, or odd if P is odd. 

Since it is difficult to decide the value of the variable as even or odd, let us consider both the options one after the other, as described earlier, in the rules of the game section.  For convenience, we always consider the even option as the first option and then come back for the computation of the odd option.   

Now, let us split the computation into two scenarios (i) P-even and (ii) P-odd.  

\subsubsection{\underline{Case: n = (4P+3), P-Even:}}
 
$5^{th}$ term = $Col^5$(n) = {$\underbrace{\frac{(9P + 8)}{2}}$} = \{{$\underbrace{(4P+4)}_{Even\ only}$} + {$\underbrace{(\frac{P}{2})
}_{Odd\ or\ Even ?} $\}} 

This $5^{th}$ term is an even number iff ($\frac{P}{2}$) is even.

\subsubsection{\underline{Case: n = (4P+3), P - Even, $(\frac{P}{2})$ - Even:}}

$6^{th}$ term = $Col^6$(n) = {$\underbrace{\frac{(9P + 8)}{4}}$} = \{$\underbrace{(2P+2) + (\frac{P}{4})}$\} $<$ n \{=(4P+3)\}

The value of $6^{th}$ term is less than the initial number (n), i.e., $Col^6$(n) = $\frac{(9n+5)}{16} < $ n.   This implies, the Collatz sequence for n=(4P+3) converges iff P is of even value and also, ($\frac{P}{2}$) is of even value\footnote{If $6^{th}$ term is the PSO term, then the condition to be considered for the evaluation of the value of the variable is $5^{th}$ term (previous term).  This is because, if this whole term value is even, then only we divide the term by 2 and the term value becomes less than the initial number.  Hence, this term provides the condition for the values of the variable to be evaluated}.  

That is to say, if P is a multiple of 4, then the sequence reaches a pre-proven sequence at the $6^{th}$ term (PSO=6).  Since the sequence converges to a value lower than the initial number (n) and the number lower than n is already proven to converge in this sequential computation method, the sequence converges to a term value of 1 with the pre-proven sequence order of PSO = 6, for n = (4P+3) with P = even, ($\frac{P}{2}$) = even, i.e., iff P is a multiple of 4. 

It should be noted here that once ($\frac{P}{2}$) is even, P will obviously be even and hence, only the last condition verification is sufficient to verify the convergence for the given number.  However, in case of any doubt, it may be a better option to verify all the conditions from the beginning till the PSO term is reached, just for a complete confirmation.

\textbf{The Collatz sequence for the number n = 4P+3, where ($\frac{P}{2}$) is even, converges with PSO = 6.}  Hence, any positive integer n = (4P+3) that corresponds to a value of P which is divisible by 4 converges. i.e., if P = 2*(2Y)=4Y, then, n = (4P+3) = (16Y+3) converges for all Y = \{0, 1, 2, 3, ....\}.

\subsubsection{\underline{Case: n = (4P+3), P - Even, $(\frac{P}{2})$ - Odd:}}

Till we use the P-even condition as valid, we could borrow the computation done earlier and use those information directly in subsequent computation.  Hence, for $(\frac{P}{2})$ - Odd, 

$5^{th}$ term = $Col^5$(n) = {$\underbrace{\frac{(9P + 8)}{2}}$} = \{{$\underbrace{(4P+4)}_{Even\ only}$} +  {$\underbrace{(\frac{P}{2})}_{considered\ as\ Odd} $} \}     $\therefore$ odd value term.  Then, 
 
$6^{th}$ term = $Col^6$(n) = \{$\underbrace{3*\frac{(9P + 8)}{2} + 1}$\} = {$\underbrace{\frac{(27P+26)}{2}}_{even-only\ term}$} = \{{$\underbrace{(13P + 13)}_{Odd}$} + {$\underbrace{(\frac{p}{2})}_{Odd}$\}} 

$7^{th}$ term = $Col^7$(n) = {$\underbrace{\frac{(27P+26)}{4}}$} = \{{$\underbrace{(6P + 6)}_{Even\ only}$}+{$\underbrace{\frac{(3P+2)}{4}}_{Even\ or\ Odd\ ?}$\}} - Even value, iff ($\frac{3P+2}{4}$) is even.  

\subsubsection{\underline{Case: n = (4P+3), P - Even, ($\frac{P}{2}$) - Odd, ($\frac{3P+2}{4}$) - Even:}}

$8^{th}$ term = $Col^8$(n) = $\underbrace{\cfrac{(27P+26)}{8}}$ = \{$\underbrace{(3P + 3) + \frac{(3P+2)}{8}}$\}  $<$ n \{= (4P+3)\}.  

i.e., the value of $8^{th}$ term is less than that of the intial number n = (4P+3).  \\Hence, PSO = 8, iff ($\frac{3P+2}{4}$) is Even.  
 
\textbf{ $\therefore$ The Collatz sequence for the number n = (4P+3), where P - even, ($\frac{P}{2}$) - odd and ($\frac{3P+2}{4}$) - even converges with PSO = 8.  Hence, the sequence for any positive integer that corresponds to a value of P such that ($\frac{3P+2}{4}$)- even converges.} 

Let us understand this scenario for the values of P.  If P should be even, and ($\frac{P}{2}$) should be odd, then, it should be of the form P=2*(2Y+1)=(4Y+2).  In such as case, ($\frac{3P+2}{4}$)=(3Y+2).  Then, Y should be even for the term ($\frac{3P+2}{4}$) to be even.  If we say, Y=2x, as Y should be even only, then P becomes =(8x+2).  Therefore, n=(4P+3), with P = (8x+2), converges with a PSO value of 8, for all x = \{0, 1, 2, 3, ....\}.  

On similar lines, if we re-write the first case where ($\frac{P}{2}$) should be even, then those numbers could be written as n = (4P+3) with P = 4Y where Y could be even or odd,  i.e., Y = 2x or Y = (2x+1).  Hence, it turns out that P = (8x+0) and (8x+4), for all x = \{0, 1, 2, 3, ....\}\footnote{As obvious, generally the  even numbers are represented as (2x) and odd numbers as (2x+1), where x is any positive integer, including zero, i.e., x = \{= 0, 1, 2, 3, ...\}.  To cover all the integers of the natural numbers of the form n = (4P+3), if we split the set of P values as octets, then the possible values of P are P = \{ (8x+0), (8x+1), (8x+2), (8x+3), (8x+4), (8x+5), (8x+6) and (8x+7) \}.  Among them, the even values of P correspond to P(Even) = \{ (8x+0), (8x+2), (8x+4) and (8x+6)\} whereas the odd values of P correspond to P(Odd) = \{ (8x+1), (8x+3), (8x+5) and (8x+7) \}.}.

Here, if P = (8x+2), ($\frac{3P+2}{4}$) = (6x+2) - even only whereas if P = (8x+6), ($\frac{3P+2}{4}$) = (6x+5) - odd only.  Hence, among the numbers satisfying the conditions of P-even, ($\frac{P}{2}$) - odd, ($\frac{3P+2}{4}$) - odd, are the numbers of the form P = (8x+6), considering P = (4Y+2) where Y = (2x+1).  Hence, we still need to explore the option of ($\frac{3P+2}{4}$) - Odd.   

With this, we have moved from the quartet analysis to octet analysis.   

\subsubsection{\underline{Case: n = (4P+3), P - Even, ($\frac{P}{2}$) - Odd, ($\frac{3P+2}{4}$) - Odd:}}

$7^{th}$ term = $Col^7$(n) = {$\underbrace{\frac{(27P+26)}{4}}$} = {$\underbrace{(6P + 6)}_{Even\ only}$} + {$\underbrace{\frac{(3P+2)}{4}}_{Even\ or\ Odd\ ?}$} :  Odd number, if ($\frac{3P+2}{4}$) is odd. 

Then, 

$8^{th}$ term = $Col^8$(n) = \{$\underbrace{3*(\frac{27P+26}{4}) +1}$\} = {$\underbrace{\frac{(81P+82)}{4}}_{Even\ only}$} = 
\{$\underbrace{(20P+20) + (\frac{P+2}{4})}$\} 

$9^{th}$ term = $Col^9$(n) = {$\underbrace{\frac{(81P+82)}{8}}$} = {$\underbrace{(10P + 10)}_{Even\ only}$} +{$\underbrace{\frac{(P+2)}{8}}_{Even\ or\ Odd\ ?}$}:  Even number, iff ($\frac{P+2}{8}$) is even.  

$10^{th}$ term = $Col^{10}$(n) = {$\underbrace{\frac{(81P+82)}{16}}$} = {$\underbrace{(5P + 5)}_{Odd\ only}$} + {$\underbrace{\frac{(P+2)}{16}}_{Even\ or\ Odd\ ?}$} :  Even number, if ($\frac{P+2}{16}$) is odd. 

$11^{th}$ term = $Col^{11}$(n) = {$\underbrace{\frac{(81P+82)}{32}}$} = {$\underbrace{(2P + 2)}$} + {$\underbrace{\frac{(17P+18)}{32}}$} $<$ n \{= (4P+3)\}. 

i.e., the value of $11^{th}$ term in the sequence is less than that of the intial number n = (4P+3), iff ($\frac{P+2}{16}$) is Odd.  Hence, PSO = 11.  

\textbf{ $\therefore$ The Collatz sequence for the number n = (4P+3), where P - even, ($\frac{P}{2}$) - odd, and ($\frac{3P+2}{4}$) - odd, ($\frac{P+2}{16}$) - odd, converges with PSO = 11.  Hence, the sequence for any positive integer that corresponds to a value of P such that ($\frac{P+2}{16}$)- odd converges.} 

Considering P = (8x+6) as the only option remaining for the given condition, as explained earlier, ($\frac{P+2}{16}$) = ($\frac{8x+8}{16}$) = ($\frac{x+1}{2}$) should be odd.   $\therefore$ the allowed values of x for ($\frac{x+1}{2}$) to be odd is x = \{1, 5, 9, 13, ....\}.  This implies, P = (8x+6), x = (4v+1), for all v = \{0, 1, 2, 3, ....\}.   

Now, the computation leads to further split of number sets and the tail would keep growing, as we progress further.  Hence, let us hold on with other P-even options and explore the situation with P-odd number options and continue the computation for n = (4P+3) with P = (8x+6) later.  

\subsubsection{\underline{Case: n = (4P+3), P-Odd:}} 

$4^{th}$ term = $Col^4$(n) = {$\underbrace{\frac{(18P + 16)}{2}}$} = {$\underbrace{(9P+8)}_{Odd\ only}$}.  

Now, considering P as odd, the fourth term $Col^4$(n) is an odd value term.  Hence, 

$5^{th}$ term = $Col^5$(n) = \{$\underbrace{3*(9P + 8) + 1}$\} = {$\underbrace{(27P+25)}$} - Even only term. 

$6^{th}$ term = $Col^6$(n) = {$\underbrace{\frac{(27P + 25)}{2}}$} = {$\underbrace{(13P + 12)}_{Odd\ term}$} + {$\underbrace{\frac{(P+1)}{2}}_{Even\ or\ Odd\ ?}$}

Now, considering the whole term as even, the part ($\frac{P+1}{2}$) should be of odd value.  

$7^{th}$ term = $Col^7$(n) = {$\underbrace{\frac{(27P+25)}{4}}$} = \{{$\underbrace{(6P + 6)}_{Even\ only}$} + {$\underbrace{\frac{(3P+1)}{4}}_{Even\ or\ Odd\ ?}$\}} -- Even, iff ($\frac{3P+1}{4}$) is even. 

$8^{th}$ term = $Col^8$(n) = {$\underbrace{\frac{(27P+26)}{8}}$} = 
 \{$\underbrace{(3P + 3) + \frac{(3P+1)}{8}}$\}, $<$ n \{= (4P+3)\}. 

i.e., the value of $8^{th}$ term is less than that of the initial number n = (4P+3).  Hence, PSO = 8, iff ($\frac{3P+1}{4}$) is Even.

\textbf{ $\therefore$ The Collatz sequence for the number n = (4P+3), where P - odd, ($\frac{P+1}{2}$) - odd, and ($\frac{3P+1}{4}$) - even, converges with PSO = 8.  Hence, the sequence for any positive integer that corresponds to a value of P such that ($\frac{3P+1}{4}$)- even converges.} 

This condition necessitates, from the octet set of numbers, that P = (8x+5). i.e., ($\frac{3P+1}{4}$) = (6x+4), which is even for all integer values of x.  Hence, the Collatz sequence for odd numbers of the form n = (4P+3) with P = (8x+5) converges for all values of x = \{0, 1, 2, 3, ...\}.  

Let us continue with the computation for n = (4P+3) for the remaining sets of P = (8x+c), considering all of the sets together with the octet combinations, in the next section\footnote{It should be noted here that when the numbers are split into various sets of the form n = (b*V + c), where V is a variable, c should have values ranging from 0 upto (b-1), i.e., for n = (4P + c), c = \{0, 1, 2, 3\} whereas for P = (8x + c), c=\{0, 1, 2, 3, 4, 5, 6, 7\} to include all the numbers without any exclusion.  This is the reason why the quartets contained n=\{(4P+0) to (4P+3)\} whereas now octets should contain P=\{(8x+0)-to-(8x+7)\}.  The split is purely for computation simplification and to manage the exploration easily.}.  

In the mean time, it is necessary to understand the process of evaluating the values of integers (P) those converge to a specific condition of even or odd combination of sets of numbers, to eliminate the confirmed-convergence and to find the remaining set of numbers to be proven in a systematic way.  

\subsection{P-value Evaluation (Example):} 

At every stage of the computation, the possible values of integer sets those are converging are to be evaluated (octet analysis) and confirmed.  This process is presented here with two examples. 

\textbf{Example I:} 

Let us consider the case of PSO = 8 if ($\frac{3P+2}{4}$) is even: 

This is the condition reached when computing the convergence for n = (4P+3) and while doing so, we have come across other conditions such as P - Even, ($\frac{P}{2}$) - Odd, and then ($\frac{3P+2}{4}$) - Even to reach PSO = 8, at steps 5, 6 and 7, respectively.  

Since we are considering the set of integers as octets, the possibilities and non-possibilities are separated from the set of all P values, viz., P = \{(8x+0), (8x+1), (8x+2), (8x+3), (8x+4), (8x+5), (8x+6) and (8x+7)\}.  

Among these sets, if we have to consider only those sets for which P are Even, then the options available are: P - Even = \{(8x+0), (8x+2), (8x+4) and (8x+6) only\}.

Among these sets, if we have to consider only those for which ($\frac{P}{2}$) are Odd, then the options available are: P = \{(8x+2) and (8x+6)\}.

In the next step, when the condition required is ($\frac{3P+2}{4}$) - Odd, then we have to find the value by substituting the value of P in the expression.

\begin{equation*} 
{ \left(\frac{3P+2}{4} \right) } = 
\begin{cases}
(6x+2) & \text{if P = (8x+2)}\\ (6x+5) & \text{if P = (8x+6)}
\end{cases}
\end{equation*} 

Hence, we conclude that only the set of numbers of the form n = (4P+3) with P = (8x+2), for all values of x = \{0, 1, 2, 3, ...\} converges to a PSO value of 8.  This is a simple and direct evaluation.  

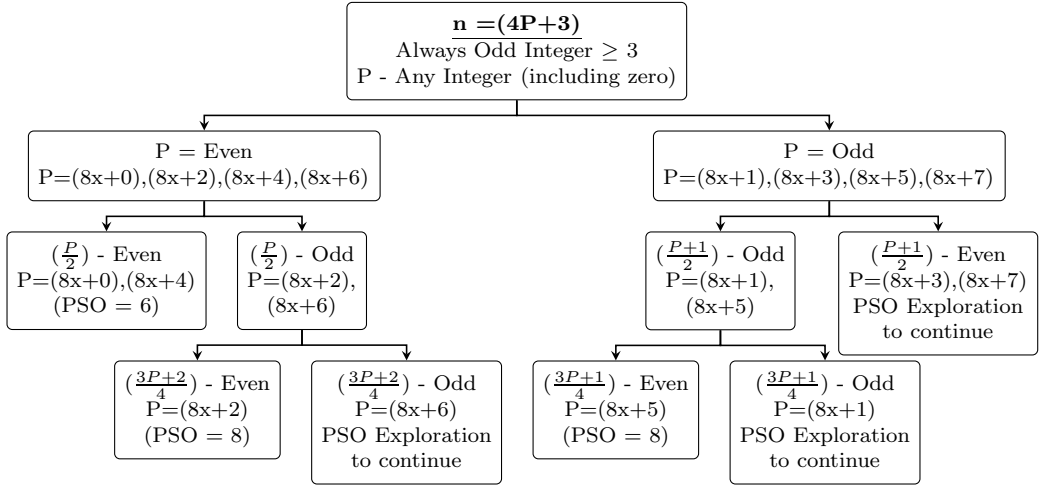
\begin{center}
\begin{forest} 
      for tree={	font=\scriptsize, align=center, 
      edge={draw, semithick, -stealth, rectangle},       
      anchor=north, forked edge,  grow=south, 
      s sep=4mm,   l sep=4mm,    fork sep=2mm, 
   tier/.option=level,  } 
   [{\textbf{\underline{n =(4P+3)}}\\ Always Odd Integer $\ge$ 3 \\P - Any Integer (including zero)}, basic
       [{P = Even \\P=(8x+0),(8x+2),(8x+4),(8x+6)}, basic 
         [{($\frac{P}{2}$) - Even \\ P=(8x+0),(8x+4) \\ (PSO = 6)}, basic ] 
         [{($\frac{P}{2}$) - Odd \\P=(8x+2),\\ (8x+6)}, basic
         [{($\frac{3P+2}{4}$) - Even\\P=(8x+2)\\(PSO = 8)}, basic]
         [{($\frac{3P+2}{4}$) - Odd \\P=(8x+6)\\ PSO Exploration \\to continue}, basic]]  ]
      [{P = Odd \\P=(8x+1),(8x+3),(8x+5),(8x+7)}, basic
         [{($\frac{P+1}{2}$) - Odd \\P=(8x+1),\\ (8x+5)}, basic
        [{($\frac{3P+1}{4}$) - Even\\P=(8x+5)\\(PSO = 8)}, basic ] 
         [{($\frac{3P+1}{4}$) - Odd \\P=(8x+1)\\ PSO Exploration \\to continue}, basic]] 
         [{($\frac{P+1}{2}$) - Even \\ P=(8x+3),(8x+7) \\PSO   Exploration \\to continue}, basic]] ] 
\end{forest}  
\end{center} 

Another way of arriving at the same result is as follows: \\ If ($\frac{3P+2}{4}$) should be even, then the permissible values of P are \{2, 10, 18, 26, ...(+8) ...\}.  This may be written in the octet form as P = (8x+2), where x = \{0, 1, 2, 3, ...\}, for all positive integer values, including zero. 

As we progress with these computations, we prefer to use this second option of computing the permissible values as the process appears to be  more simpler to manage.  

\textbf{Example II:}

Next, let us consider an example from a complex scenario.  

When we compute the sequence for n = (4P+3), P = (8x+1), we arrive at the condition, ($\frac{9x+11}{16}$)-even.  Then, the verification is done by exploring all possible values of x for which this condition is true.  The permissible values of x in this case are x = \{13, 29, 45, ... (+16)...\}.  
$\therefore$ If ($\frac{9x+11}{16}$) = even, then, we could write the same in the form of next octet as x = (8q+5), wherein q = \{all odd numbers\} = \{1, 3, 5, ... \}.  

Because of this evaluation step in the process of confirmation of the set of numbers those are converging under specific conditions, the overall computation is lengthy and time consuming.   As explained till now, 87.50\% of all the positive integers have been confirmed to converge to the PSO terms and hence, the convergence to 1.  However, the remaining 12.50\% , that too the last $<$ 1\% of the proof might take a long time to complete due to the laborious nature of computation\footnote{By this procedure, presently, $>$ 99\% integers have been found to converge to their PSO term values and only $<$ 1\% of all the natural numbers (positive integers) are remaining to be confirmed.  Due to the computational limitation, procedure implementation is consuming time.}     

At this stage, we could say that \textbf{the only procedure that could help in confirming the non-existence of counter-example has been established}.

\section{Evaluation Leads to Elimination}     

From the quartet analysis, it may be observed that the remaining group of numbers to be proven for convergence of Collatz sequences is of the form n=(4P+3).  Even among them, half of them in the octet are converging quickly with low PSO values, as seen in the previous section.  Therefore, the remaining set of numbers those need to be proved are n=(4P+3) where P=\{(8x+1), (8x+3), (8x+6) and (8x+7)\}.  

This octet grouping is described in the following tree-like structure.  
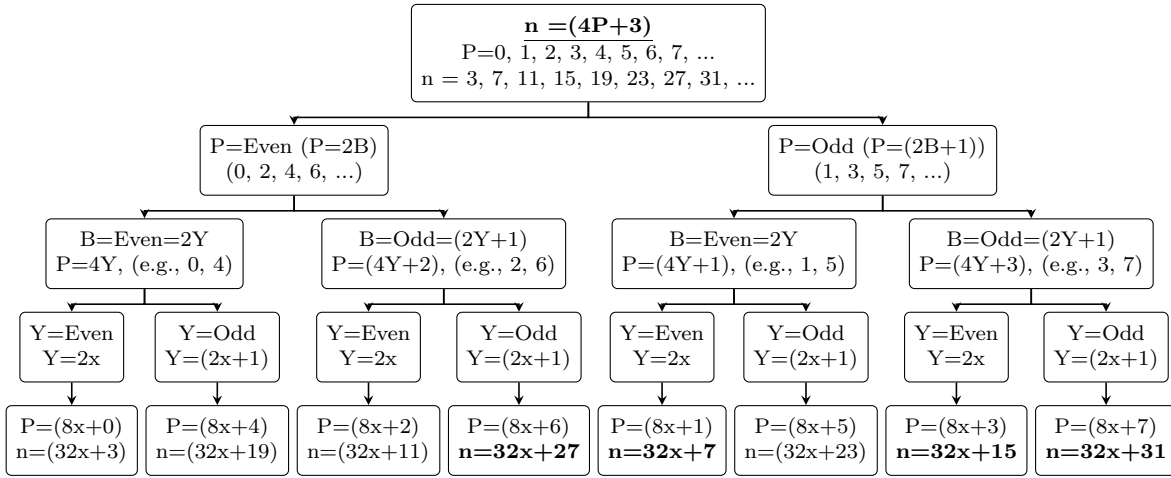
\begin{center}
\begin{forest} 
      for tree={	font=\scriptsize, align=center,  
      edge={draw, semithick, -stealth, rectangle},       
      anchor=north, forked edge,  grow=south, 
      s sep=1mm,   l sep=3mm,    fork sep=2mm, 
   tier/.option=level,  }  
   [{\textbf{\underline{n =(4P+3)}}\\P=0, 1, 2, 3, 4, 5, 6, 7, ... \\n = 3, 7, 11, 15, 19, 23, 27, 31, ...}, basic 
     [{P=Even (P=2B)  \\(0, 2, 4, 6, ...)}, basic
         [{B=Even=2Y \\P=4Y, (e.g., 0, 4)}, basic 
             [{Y=Even \\Y=2x}, basic 
                 [{P=(8x+0) \\n=(32x+3)}, basic]     ]
             [{Y=Odd \\Y=(2x+1)}, basic 
                 [{P=(8x+4) \\n=(32x+19)}, basic] ] 
         ] 
         [{B=Odd=(2Y+1) \\P=(4Y+2), (e.g., 2, 6)}, basic
            [{Y=Even \\Y=2x}, basic 
                [{P=(8x+2) \\n=(32x+11)}, basic] ]
            [{Y=Odd \\Y=(2x+1)}, basic 
                [{P=(8x+6) \\\textbf{n=32x+27}}, basic] ] 
         ]
      ] 
      [{P=Odd (P=(2B+1)) \\(1, 3, 5, 7, ...)}, basic
         [{B=Even=2Y \\P=(4Y+1), (e.g., 1, 5)}, basic
             [{Y=Even \\Y=2x }, basic 
                 [{P=(8x+1) \\\textbf{n=32x+7}}, basic] ]   
             [{Y=Odd \\Y=(2x+1)}, basic 
                 [{P=(8x+5) \\n=(32x+23)}, basic] ] 
         ] 
         [{B=Odd=(2Y+1) \\P=(4Y+3), (e.g., 3, 7)}, basic
             [{Y=Even \\Y=2x}, basic 
                 [{P=(8x+3) \\\textbf{n=32x+15}}, basic]  ] 
             [{Y=Odd \\Y=(2x+1)}, basic 
                 [{P=(8x+7) \\\textbf{n=32x+31}}, basic] ] ] ] ] 
\end{forest}
\end{center} 

Here onwards, each individual set of P values will be analyzed for the convergence of the Collatz sequences for those individual set of numbers.  

\subsection{Alternate Option: x as variable} 

While we computed the convergence of the set of numbers of the form n = (4P+3), we found the convergence for the numbers of the form P = (8x+0), (8x+4), (8x+2) and (8x+5), using P as variable.  But, the computation becomes simpler with the direct computation with x as variable.  However, it was necessary during the initial stages to do those computations to realize the availability of this alternate options.  Now, let us repeat one of those sets P=(8x+5) of computations using x as variable and compare the same with P as variable and demonstrate the equivalence and convenience, to adopt this alternate option for further computations.  Since the coefficient of P gets multiplied by a factor of 8 while using x as vaiable, the computation here gets simplified by a factor of $2^3$ = 8 and hence, there is no fractional part remaining in the terms of the sequence, at the initial level.  
 
\begin{center}
\begin{tabular}{||c||c|c||c||c|c||} 
\hline 
\multicolumn{6}{||c||}{Table 4: PSO = 8 Scenario Computation with both x- as well as P- as variables} \\  
\multicolumn{6}{||c||}{  } \\ 
\multicolumn{6}{||c||}{Here, P - Odd, ($\frac{P+1}{2}$) - Odd and ($\frac{3P+1}{4}$) - Even : possibilities are considered.} \\  
\hline 
\multicolumn{1}{||c||}{ } 
&\multicolumn{2}{c||}{n=(4P+3); P=(8x+5)} 
&\multicolumn{1}{c||}{Term} 
&\multicolumn{2}{c||}{n = (4P+3)}  \\ 
\hline 
\multicolumn{1}{||c||}{$n_a$} 
  &\multicolumn{1}{c|}{Formula} 
  &\multicolumn{1}{c||}{Term}  
  &\multicolumn{1}{c||}{Type}  
  &\multicolumn{1}{c|}{Formula} 
  &\multicolumn{1}{c||}{Term} \\ 
\hline
$n_0$ & n = (32x + 23) & (32x+23) & Odd & n = (4P + 3) & (4P + 3) \\ 
\hline 
$n_1$ & \{3*(32x+23)+1\}  & (96x+70) & Even & \{3*(4P+3)+1\}  & (12P + 10) \\ 
$n_2$ & (96x + 70)$\div$ 2 & (48x+35) & Odd & (12P + 10) $\div$ 2 & (6P + 5) \\ 
$n_3$ & \{3*(48x+35)+1\} & (144x+106) & Even & \{3*(6P+5)+1\}  &  (18P + 16) \\ 
$n_4$ & (144x + 106) $\div$ 2 & (72x+53) & Odd & (18P + 16) $\div$ 2 & (9P + 8) \\ 
$n_5$ & \{3*(72x+53)+1\} & (216x+160) & Even & \{3*(9P+8)+1\} &  (27P + 25) \\ 
$n_6$ & (216x+160) $\div$ 2  & (108x+80) & Even & (27P + 25) $\div$ 2 & (13P+12)+($\frac{P+1}{2}$) \\ 
$n_7$ & (108x + 80) $\div$ 2 & (54x+40) & Even & (27P + 25) $\div$ 4 & (6P+6)+($\frac{3P+1}{4}$) \\ 
$n_8$ & (54x + 40) $\div$ 2 & (27x+20) & \{$<$n\} & (27P + 25) $\div$ 8 & (3P+3)+($\frac{3P+1}{8}$) \\ 
\hline
\end{tabular}  
\end{center}  

\textbf{It should be noted that in these computations all the variables (n, P or x) are all integers and the value of the terms in the sequences are also integers only and no fractions.}  While using variables, the term is split into integer and fraction portions whereas once the value of the variable is substituted, the term value will only be integers and no fraction.

By comparison, both the methods lead to the same result. Hence, it is the choice we could utilize to compute the remaining part of the analysis.  As it is only the substitution of the value of P as a function of x, the comparison is obvious.  Hence, the computation of proven sequence reach term could be carried out by different starting point.  These are, for example,  

\begin{itemize} 

\item{n = (4P+3) with P as variable. } 

\item{n = (8B+3) and n = (8B+7), with B as variable, where P = 2B (even) or P = (2B+1) (odd). }

\item{n = (4P+3) where P = (8x+0), (8x+2), (8x+4) and (8x+6) - even values and P = (8x+1), (8x+3), (8x+5) and (8x+7) - odd values, i.e.,  n = (32x+3), (32x+7), (32x+11), (32x+15), (32x+19), (32x+23), (32x+27) and (32x+31) respectively, with x as variable.} 

\end{itemize} 

All the procedures will lead to the same results and hence, the choice is considered by the convenience of the person computing the sequence.  In the present computations, we have already confirmed the convergence of the set of numbers of the form n = (4P+3) where P = (8x+0), (8x+2) and (8x+4) - even values and P = (8x+5) - odd values.  Now, for the remaining set of numbers of the form n = (4P+3) where P = (8x+6) - even values, and P = (8x+1), (8x+3) and (8x+7) - odd values, i.e., n = (32x+7), (32x+15), (32x+27) and (32x+31), let us progress with x as variable, here onwards. 

\subsection{Computation of Sequences with x - as variable} 

At this stage, the remaining 12.5\% of the total positive integers are the ones those are posing complex scenarios wherein the proof of the conjecture takes a long time, due to laborious computation.  As we could see further, the last $<$ 3\% of the numbers are the ones needing special computation facility to confirm the convergence of the type of set of numbers those are converging with various PSO values. 

Since all the computational procedures have been established and explained with examples, details of how to compute are available in the previous sections.  Now onwards, the individual steps of computation and the convergence scenarios would be tabulated for simplification.  Finally, the summary of results to take stock of the progress would be presented at the end. 

\subsubsection{System: n=(32x+27) [i.e., n=(4P+3); P=(8x+6)]} 

Among all the positive \textbf{\underline{even values of P}}, the  remaining set to be evaluated is P = (8x+6); i.e., n=(32x+27), where x =\{0, 1, 2, 3, ...(integers including zero)\}.   

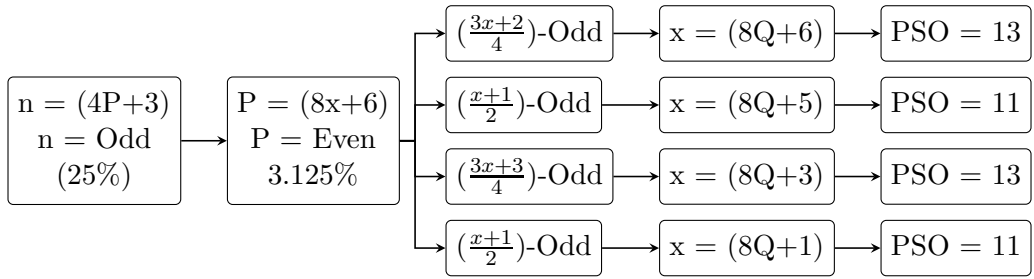
\begin{center}
\begin{tabular}{||c|c|c|c|c||} 
\hline 
\multicolumn{5}{||c||}{\textbf{Table 5.1}: n=(4P+3); P=(8x+6): Computation with x - as variable} \\  
\hline 
\multicolumn{5}{||c||}{\textbf{n = (32x+27); x - Odd}} \\
\hline 
\multicolumn{1}{||c|}{$n_a$} 
  &\multicolumn{1}{|c|}{Formula} 
  &\multicolumn{1}{|c|}{Term} 
  &\multicolumn{1}{|c|}{Type} 
  & \multicolumn{1}{|c||}{Condition} \\ 
\hline
$n_1$ & \{3*(32x+27)+1\}      & (96x + 82)   & Even  & - \\ 
$n_2$ & (96x + 82)$\div$ 2    & (48x + 41)   & Odd   & - \\ 
$n_3$ & \{3*(48x+41)+1\}      & (144x + 124) & Even  & - \\ 
$n_4$ & (144x + 124) $\div$ 2 & (72x + 62)   & Even  & - \\ 
$n_5$ & (72x + 62)$\div$ 2    & (36x + 31)   & Odd   & - \\ 
$n_6$ & \{3*(36x+31)+1\}      & (108x + 94)  & Even  & - \\ 
$n_7$ & (108x + 94) $\div$ 2  & (54x + 47)   & Odd   & - \\ 
$n_8$ & \{3*(54x + 47)+1\}    & (162x + 142) & Even  & - \\ 
$n_9$ & (162x + 142) $\div$ 2  & (81x + 71)  & Even  & if x is Odd \\
$n_{10}$ & (81x + 71) $\div$ 2  & (40x + 35) + ($\frac{x+1}{2}$) & Even & ($\frac{x+1}{2}$)-Odd \\
$n_{11}$ & (81x + 71) $\div$ 4  & (20x + 17) + ($\frac{x+3}{4}$) & \{$<$ n\} & ($\frac{x+1}{2}$)-Odd \\
\hline 
\multicolumn{5}{||c||}{Result: Since $n_{11} < n $, \textbf{PSO = 11}, if \{x - Odd\} and \{($\frac{x+1}{2}$) - Odd\}} \\ 
\multicolumn{5}{||c||}{\textbf{i.e., x = (8q+1) or (8q+5), for all q = \{0, 1, 2, 3, ... .\}.}}\\ 
\hline 
\multicolumn{5}{||c||}{   } \\ 
\multicolumn{5}{||c||}{\textbf{Table 5.2: n = (32x+27); x - Odd}} \\  
\hline 
\multicolumn{1}{||c|}{$n_a$}  
  &\multicolumn{1}{|c|}{Formula} 
  &\multicolumn{1}{|c|}{Term} 
  &\multicolumn{1}{|c|}{Type} 
  & \multicolumn{1}{|c||}{Condition} \\ 
\hline 
$n_9$ & (162x + 142) $\div$ 2   & (81x + 71)   & Even   &  ($\because$ x - Odd) \\
\hline 
$n_{10}$ & (81x + 71) $\div$ 2 & (40x + 35) + ($\frac{x+1}{2}$) & Odd & ($\frac{x+1}{2}$) - Even \\
$n_{11}$ & \{3*($\frac{81x + 71}{2}$)+1\} & (121x + 107) + ($\frac{x+1}{2}$) & Even & ($\frac{x+1}{2}$) - Even \\
$n_{12}$ & (243x + 215) $\div$ 4 & (60x + 53) + ($\frac{3x+3}{4}$) & Even & ($\frac{3x+3}{4}$) - Odd \\ 
$n_{13}$ & (243x + 215) $\div$ 8 & (30x + 26) + ($\frac{3x+7}{8}$) & $< n$ & ($\frac{3x+3}{4}$) - Odd \\
\hline
\multicolumn{5}{||c||}{Result: Since $n_{13} < n $, \textbf{PSO = 13}, if \{($\frac{3x+3}{4}$) - Odd\}} \\ 
\multicolumn{5}{||c||}{\textbf{i.e., x = (8q+3), for all q = \{0, 1, 2, 3, ...\}.}}\\ 
\hline 
\multicolumn{5}{||c||}{  } \\
\multicolumn{5}{||c||}{\textbf{Table 5.3: n = (32x+27); x - Even}} \\
\hline 
\multicolumn{1}{||c|}{$n_a$} 
  &\multicolumn{1}{|c|}{Formula} 
  &\multicolumn{1}{|c|}{Term} 
  &\multicolumn{1}{|c|}{Type} 
  & \multicolumn{1}{|c||}{Condition} \\ 
\hline
$n_9$ & (162x + 142) $\div$ 2 & (81x + 71)   & Odd   &  ($\because$ x - Even) \\
\hline 
$n_{10}$ & \{3*(81x + 71)+1\} & (243x + 214) & Even  &  Even only\\
$n_{11}$ & (243x + 214) $\div$ 2  & (121x + 107) + ($\frac{x}{2}$) & Even & if ($\frac{x}{2}$) - Odd \\
$n_{12}$ & (243x + 214) $\div$ 4  & (60x + 53) + ($\frac{3x+2}{4}$) & Even & ($\frac{3x+2}{4}$) - Odd \\
$n_{13}$ & (243x + 214) $\div$ 8  & (30x + 26) + ($\frac{3x+6}{8}$) & \{$<$ n\} & ($\frac{3x+2}{4}$) - Odd \\
\hline
\multicolumn{5}{||c||}{Result: Since $n_{13} < n$, \textbf{PSO = 13}, if \{($\frac{3x+2}{4}$) - Odd\}} \\ 
\multicolumn{5}{||c||}{\textbf{i.e., x = (8q+6), for all q = \{0, 1, 2, 3, ...\} converges.}}\\  
\hline 

\end{tabular} 
\end{center}  

\begin{center}
\begin{forest} 
      for tree={	align=center,  
      	edge={draw, semithick, -stealth, rectangle},       
      	anchor=west,  forked edge,  grow=east, 
      	s sep=2mm,         	l sep=6mm,   fork sep=2mm,   
      	tier/.option=level, }  
       [{n = (4P+3)\\n = Odd \\(25\%)}, basic
          [{P = (8x+6) \\P = Even \\3.125\%}, basic 
             [{($\frac{x+1}{2}$)-Odd}, basic [{x = (8Q+1)}, basic [{PSO = 11}, basic] ] ]
             [{($\frac{3x+3}{4}$)-Odd}, basic [{x = (8Q+3)}, basic [{PSO = 13}, basic] ] ]
             [{($\frac{x+1}{2}$)-Odd}, basic [{x = (8Q+5)}, basic [{PSO = 11}, basic] ] ]
             [{($\frac{3x+2}{4}$)-Odd}, basic [{x = (8Q+6)}, basic [{PSO = 13}, basic] ]]] ]
\end{forest}
\end{center} 

With this, it is clear that for n=(4P+3), 50\% of the total numbers in the set P = (8x+6) converge to a PSO value of either 11 or 13.  The remaining set of numbers within the set of numbers of the form P = (8x+6) are \{x = (8Q+0), x = (8Q+2), x = (8Q+4), x = (8Q+7)\}. 
 
Among all the positive \textbf{\underline{odd values of P}}, the  remaining sets to be evaluated are P = (8x+1), (8x+3) and (8x+7), where x =\{0, 1, 2, 3, ...(integers including zero)\}.   

\subsubsection{System: n=(32x+7); [i.e., n=(4P+3); P=(8x+1)] } 

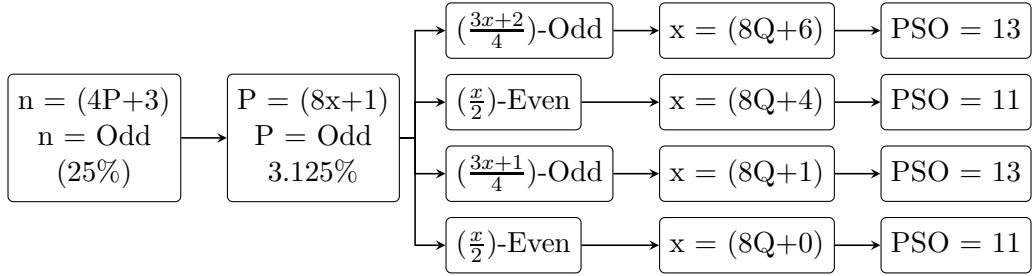
\begin{center} 
\begin{forest} 
      for tree={	align=center,  
      	edge={draw, semithick, -stealth, rectangle},       
      	anchor=west,  forked edge,  grow=east, 
      	s sep=2mm,         	l sep=6mm,   fork sep=2mm,   
      	tier/.option=level, }  
       [{n = (4P+3)\\n = Odd \\(25\%)}, basic
          [{P = (8x+1) \\P = Odd \\3.125\%}, basic 
             [{($\frac{x}{2}$)-Even}, basic [{x = (8Q+0)}, basic [{PSO = 11}, basic] ] ]
             [{($\frac{3x+1}{4}$)-Odd}, basic [{x = (8Q+1)}, basic [{PSO = 13}, basic] ] ]
             [{($\frac{x}{2}$)-Even}, basic [{x = (8Q+4)}, basic [{PSO = 11}, basic] ] ]
             [{($\frac{3x+2}{4}$)-Odd}, basic [{x = (8Q+6)}, basic [{PSO = 13}, basic] ]]] ]
\end{forest}
\end{center} 

\subsubsection{System: n=(32x+15); [i.e., n=(4P+3); P=(8x+3)] } 

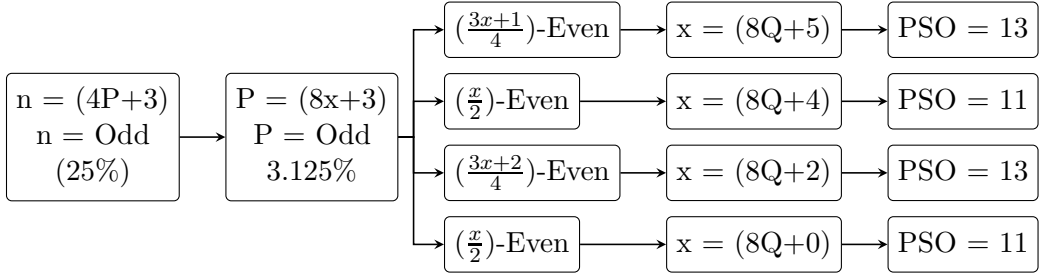
\begin{center} 
\begin{forest} 
      for tree={	align=center,  
      	edge={draw, semithick, -stealth, rectangle},       
      	anchor=west,  forked edge,  grow=east, 
      	s sep=2mm,         	l sep=6mm,   fork sep=2mm,   
      	tier/.option=level, }  
       [{n = (4P+3)\\n = Odd \\(25\%)}, basic
          [{P = (8x+3) \\P = Odd \\3.125\%}, basic 
             [{($\frac{x}{2}$)-Even}, basic [{x = (8Q+0)}, basic [{PSO = 11}, basic] ] ]
             [{($\frac{3x+2}{4}$)-Even}, basic [{x = (8Q+2)}, basic [{PSO = 13}, basic] ] ]
             [{($\frac{x}{2}$)-Even}, basic [{x = (8Q+4)}, basic [{PSO = 11}, basic] ] ]
             [{($\frac{3x+1}{4}$)-Even}, basic [{x = (8Q+5)}, basic [{PSO = 13}, basic] ]]] ]
\end{forest}
\end{center} 

\subsubsection{System: n=(32x+31); [i.e., n=(4P+3); P=(8x+7)] } 

\begin{center}
\begin{forest} 
      for tree={	align=center,  
      	edge={draw, semithick, -stealth, rectangle},       
      	anchor=west,  forked edge,  grow=east, 
      	s sep=2mm,         	l sep=6mm,   fork sep=2mm,   
      	tier/.option=level, }  
       [{n = (4P+3)\\n = Odd \\(25\%)}, basic
          [{P = (8x+7) \\P = Odd \\3.125\%}, basic 
             [{($\frac{3x+2}{4}$)-Even}, basic [{x = (8Q+2)}, basic [{PSO = 13}, basic] ] ] ] ]
\end{forest}
\end{center} 

The following compilation gives the status upto PSO=32.  

\begin{center} 
\begin{tabular}{||c||c|c||c|c||c|c||c|c||} 
\multicolumn{9}{c}{\textbf{Table 6: The PSO counts at various PSO stages.}} \\[5pt] 
\hline 
\multicolumn{1}{||c||}{PSO} 
&\multicolumn{2}{c||}{P=(8x+1)} 
&\multicolumn{2}{c||}{P=(8x+3)} 
&\multicolumn{2}{c||}{P=(8x+6)} 
&\multicolumn{2}{c||}{P=(8x+7)} \\ 
\hline 
\multicolumn{1}{||r||}{n = } 
&\multicolumn{2}{c||}{n = (32x+7)} 
&\multicolumn{2}{c||}{n = (32x+15)} 
&\multicolumn{2}{c||}{n = (32x+27)} 
&\multicolumn{2}{c||}{n = (32x+31)} \\  
\hline 
x-Type & Even & Odd & Even & Odd & Even & Odd & Even & Odd \\
\hline 
11 & 1 & 0 & 1 & 0 & 0 & 1 & 0 & 0 \\ 
\hline 
13 & 1 & 1 & 1 & 1 & 1 & 1 & 1 & 0 \\ 
\hline 
16 & 1 & 2 & 1 & 2 & 2 & 1 & 2 & 1 \\ 
\hline 
19 & 2 & 5 & 2 & 5 & 5 & 2 & 5 & 4 \\ 
\hline 
21 & 5 & 14 & 5 & 14 & 14 & 5 & 14 & 14 \\ 
\hline 
24 & 9 & 28 & 9 & 28 & 28 & 9 & 28 & 34 \\ 
\hline 
26 & 23 & 76 & 23 & 76 & 76 & 23 & 76 & 3 \\ 
\hline 
29 & 43 & 151 & 43 & 151 & 151 & 43 & 151 & 228 \\ 
\hline 
32 & 113 & 412 & 113 & 412 & 412 & 113 & 412 & 665 \\ 
\hline 
\end{tabular} 
\end{center}  

Here again, we could see the pattern in terms of the PSO count values. 

The following systems are similar:

1. (P=(8x+1), x-Even) is similar to (P=(8x+3), x-Even) and to (P=(8x+6), x-Odd). 

2. (P=(8x+1), x-Odd) is similar to (P=(8x+3), x-Odd) and to (P=(8x+6), x-Even) and to (P=(8x+7), x-Even). 

Hence, the further progress of work would be as follows: 

1. P=(8x+1): \{x-Even and x-Odd\} as well as P=(8x+7): \{x-Odd\}. 
 
2. P=(8x+3) and P=(8x+6), both x-Even as well as x-Odd and also, P=(8x+7), x-Even are to be confirmed by similarity. 

As of now, the extent of confirmation completed is about 99\% of \{N\} and the pending work is about 1\%.  This status may be summarized in the following table. 

\begin{center} 
\begin{longtable}{||c|c|c|c||}
\multicolumn{4}{c}{\textbf{Table 7: The compilation of the extent of confirmation of convergence to PSO terms.}} \\[5pt] 
\hline 
\multicolumn{1}{||c|}{SN} 
&\multicolumn{1}{|c|}{Number Sets}  
&\multicolumn{1}{|c|}{\% of \{N\} completed} 
&\multicolumn{1}{|c||}{Remarks} \\  
\hline 
\endfirsthead  
\multicolumn{4}{c}{\textbf{Table 7: The compilation of the extent of confirmation of convergence} (\textit{Continued})}\\[5pt] 
\hline 
\multicolumn{1}{||c|}{PSO} 
&\multicolumn{1}{|c|}{Condition} 
&\multicolumn{1}{|c|}{Number sets} 
&\multicolumn{1}{|c||}{\% N completed} \\ 
\hline 
\endhead 
\hline 
\multicolumn{4}{r}{\small\itshape{(continued on the next page)}}\\
\endfoot 
\hline \multicolumn{4}{r}{\small \itshape{(This is just the current status compilation; work-in-progress).}} 
\endlastfoot 
1 & n= \{1, 2 and 3\}    &  -     & (Initial values) \\ 
\hline 
2 & n=(4P+0), n-Even & 25\% ( = 100/4) & (PSO = 1) \\ 
\hline 
3 & n=(4P+1), n-Odd  & 25\% ( = 100/4) & (PSO = 3) \\ 
\hline 
4 & n=(4P+2), n-Even & 25\% ( = 100/4) & (PSO = 1) \\ 
\hline 
  & \underline{n=(4P+3), n-Odd} &        &          \\ 
  & P=(8x+0), P-Even & 3.125\% ( = 25/8) & (PSO=6) \\
5 & P=(8x+2), P-Even & 3.125\% ( = 25/8) & (PSO=8) \\
  & P=(8x+4), P-Even & 3.125\% ( = 25/8) & (PSO=6) \\
  & P=(8x+5), P-Odd  & 3.125\% ( = 25/8) & (PSO=8) \\ 
\hline 
  &                      &         & Confirmation by \\
6 & P=(8x+3), P-Odd  & 3.125\% ( = 25/8) & similarity with P=(8x+1)\\
\hline 
  &                      &         & Confirmation by cross \\
7 & P=(8x+6), P-Even & 3.125\% ( = 25/8)& similarity with P=(8x+1) \\
\hline 
  & P=(8x+7), P-Odd  &          & Confirmation by similarity\\ 
8 &             x-Even & 1.5625\% ( = 25/8)*(4/8)& with P=(8x+1), x-Odd \\
\hline 
  & \underline{P=(8x+1), P-Odd} &          &   \\ 
  & x-Even, x=(8q+0) &          & (PSO=11) \\ 
9 & x-Even, x=(8q+4) & 1.5625\% ( = 25/8)*4/8) & (PSO=11) \\ 
  & x-Even, x=(8q+6) &          & (PSO=13) \\ 
  & x-Odd   x=(8q+1) &          & (PSO=13) \\ 
\hline 
   &                 &  0.146484375  & =3*(25/$2^{8}$); (PSO=16)\\ 
   &                 &  0.1708984375 & =7*(25/$2^{10}$); (PSO=19)\\
   & P=(8x+1), P-Odd &  0.23193359375 & =19*(25/$2^{11}$); (PSO=21)\\ 
   &                 &  0.1129150390625 & =37*(25/$2^{13}$);(PSO=24)\\ 
   & \{x-Even=(8q+2)\} &  0.15106201171875 & =99*(25/$2^{14}$); (PSO=26) \\  
   &        &  0.074005126953125 & =194*(25/$2^{16}$); (PSO=29) \\ 
10 &  and   &  0.050067901611328125 & =525*(25/$2^{18}$); (PSO=32)\\ 
   &        &  0.074863433837890625 & =1570*(25/$2^{19}$); (PSO=34)\\ 
   & \{x-Odd=  &  0.040376186370849609375 & =3387*(25/$2^{21}$); (PSO=37)\\ 
   &  (8q+3,5,7)\} &  0.057756900787353515625 & =9690*(25/$2^{22}$); (PSO=39)\\    
   &        &  0.030432641506195068359375 & =20423*(25/$2^{24}$); (PSO=42)\\ 
   &                      &  0.031845271587371826171875\# & =42742*(25/$2^{25}$); (PSO=44)\\  
   &                      &  0.004414655268192291259765625\# & =23701*(25/$2^{27}$); (PSO=47)\\ 
   &                      &  0.0008196569979190826416015625\# & =17602*(25/$2^{29}$); (PSO=50)\\ 
\hline  
   &                      &  0.048828125 & =1*(25/$2^{8}$); (PSO=16)\\
   &                      &  0.097656250 & =4*(25/$2^{10}$); (PSO=19)\\
   &  P=(8x+7), P-Odd &  0.1708984375 & =14*(25/$2^{11}$); (PSO=21)\\ 
11 &                      &  0.103759765625 & =34*(25/$2^{13}$); (PSO=24)\\ 
   &                      &  0.15716552734375 & =103*(25/$2^{14}$); (PSO=26)\\  
   & \{x-Odd=             & 0.08697509765625 & =228*(25/$2^{16}$);  (PSO=29)\\
   &  (8q+1,3,5,7)\}      &  0.063419342041015625 & =665*(25/$2^{18}$); (PSO=32)\\  
   &                      &  0.099945068359375 & =2096*(25/$2^{19}$); (PSO=34)\\ 
   &                      &  0.057065486907958984375 & =4787*(25/$2^{21}$); (PSO=37)\\ 
   &                      &  0.0848710536956787109375 & =14239*(25/$2^{22}$); (PSO=39)\\ 
   &                      &  0.04668533802032470703125 & =31330*(25/$2^{24}$); (PSO=42)\\
   &                      &  0.00484287738800048828125\#  & =6500*(25/$2^{25}$); (PSO=44)\\
\hline 
\multicolumn{4}{||c||}{\# still counting.  Computation is in progress. $< 1\%$ pending.} \\
\hline 

\end{longtable} 
\end{center}  

\section{Data Compilation}

The complete computation is very lengthy and hence, the data would be compiled and made available as depository.  A summary table would be provided with the PSO count values and then the complete computation tables would be available in the depository.  

\begin{center} 
\begin{tabular}{||c||c|c||c|c||c|c||} 
\multicolumn{7}{c}{\textbf{Table 8: The PSO counts at various PSO stages.}} \\[5pt] 
\hline 
\multicolumn{1}{||r||}{ }  
&\multicolumn{4}{c||}{n = (32x+7)} 
&\multicolumn{2}{c||}{n = (32x+31)} \\  
\hline 
\multicolumn{1}{||c||}{ PSO } 
&\multicolumn{2}{c||}{P=(8x+1), x - Even} 
&\multicolumn{2}{c||}{P=(8x+1), x - Odd } 
&\multicolumn{2}{c||}{P=(8x+7), x - Odd } \\  
\hline 
Value & Count & (Table 1.2.-) & Count & (Table 1.1.-) & Count & (Table 7.1.-)\\
\hline 
11 & 1 & 1              & - &   -               & - & - \\ 
\hline 
13 & 1 & 2              & 1 & 1                 & - & - \\ 
\hline 
16 & 1 & 3              & 2 & 2 - 3             & 1 & 1 \\ 
\hline 
19 & 2 & 4 - 5          & 5 & 4 - 8             & 4 & 2 - 5 \\ 
\hline 
21 & 5 & 6 - 10         & 14 & 9 - 22           & 14 & 6 - 19 \\ 
\hline 
24 & 9 & 11 - 19        & 28 & 23 - 50          & 34 & 20 - 53 \\ 
\hline 
26 & 23 & 20 - 42       & 76 & 51 - 126         & 103 & 54 - 156 \\  
\hline 
29 & 43 & 43 - 85       & 151 & 127 - 277       & 228 & 157 - 384 \\ 
\hline 
32 & 113 & 86 - 198     & 412  & 278 - 689      & 665 & 385 - 1049 \\ 
\hline  
34 & 331 & 199 - 529    & 1239 & 690 - 1928     & 2096 & 1050 - 3145\\ 
\hline 
37 & 698 & 530 - 1227   & 2689 & 1929 - 4617    & 4787 & 3146 - 7932\\ 
\hline 
39 & 1966 & 1228 - 3193 & 7724 & 4618 - 12341   & 14239 & 7933 - 22171 \\ 
\hline 
42 &  4072 & 3194 - 7265 & 16351 & 12342 - 28692 & 31330 & 22172 - 53501\\ 
\hline 
44 & 11433 & 7266 - 18698  & 31309* & 28693-60001 & 6500* & 53502 - 60001 \\  
\hline 
47 & 23701* & 18699 - 42399 & - & - & - & - \\ 
\hline 
50 & 17602* & 42400 - 60001 & - & - & - & - \\
\hline 
Total: & 60001 & -  & 60001 & - & 60001 & - \\ 
\hline 
\multicolumn{7}{||c||}{*Still computing; The pending confirmation is : 0.925012... \%.} \\
\hline 
\end{tabular}  
\end{center}

\section{\underline{Summary and Conclusions}} 

A method to prove the Collatz Conjecture is explored and the procedure established to confirm the convergence of the Collatz sequences, using a variable.  The method follows the confirmation of the convergence of the Collatz sequence for all the natural numbers as sets of numbers.  

Initially, some of the relevant definitions are provided.  Then, the remarkable properties of the Collatz sequences are analyzed and explained the reasons for the same.  

Then, the Collatz sequences for the first few natural numbers (1-200) were numerically computed, in a sequentially forward manner (ascending order of n).  Hence, all the numbers less than n have been proved to converge before beginning to prove the convergence of the sequence for the number n.  The pattern analysis of the trajectories of the sequences of these (1-200) natural numbers, generated by direct numerical computation, led to the concept of pre-proven sequence.  Then, the origin of quartets and subsequently, the octet patterns are arrived at.  The ``Even-Odd" combinations and the concept of (8*v + c) gave a clear picture for the overall process of evaluation of the convergence. 

Subsequently, individual set of numbers are analyzed and the convergence confirmations are arrived at.  The process of computing involves the concept of combing to avoid leaving any number sets unattended.  Once the computation and confirmation is complete, all the data would be compiled and the proof would be confirmed.  As the computation is laborious, the process would take considerable amount of time.  \emph{While the computation is in progress, the current status is being shared / published, in the mean time, to receive any suggestions and/or review comments about the method by the experts.}   

As the computation progresses, the percent remaining numbers to be proved /confirmed for convergence, keeps reducing.  Hence, it appears that the height of the number-set-tree to be finite and the proof would be achievable, though the computation is laborious and taking time.  It could be better and useful if we could lay hands on the data generated in the distributed computing project (2017) to see if there is any specific upper limit for the PSO value (i.e., the maximum height of the PSO tree).       

It is a great feeling to progress with the confirmation of the Collatz  convergence of the sequences for all the natural numbers and hoping to achieve the confirmation of convergence for all the numbers without any exception.  Once a suitable computation facility is accessible for the computation of the last about 1\% of the number sets, the work would progress fast and the full proof (all the data tables) would be made available.  Finally, the compilation of all the sets and the completion of the proof will be presented. 
 
The author acknowledges the support of his family, with patience, during the progress of the project to explore \& establish the method, and also while implementing the same.

\vspace{36pt} 

\textbf{About the author:} RAMACHANDRA BHAT has been working in the industrial R \& D for over 25 years in India.  His experience encompasses the basic research in Chemistry and materials science followed by product development in the category of Personal Care Cosmetics.  He was associated with Hindustan Lever Ltd. (now, Hindustan Unilever Ltd.), Balsara Home Products Ltd., Dabur Research Foundation (Dabur India Ltd.), Marico Ltd., Kemwell Biopharma Pvt. Ltd.  He studied pre-university with science (PCMB) from Sarada Vilas College, Mysore (1980) and graduation with Mathematics, Physics and Chemistry, from Yuvaraja's College, University of Mysore (1983). Then, he did his M.Sc. (1983-85) and Ph.D. (Physical Chemistry) from Indian Institute of Technology, Bombay (1986-1990).  Visited University of Warwick, Coventry, UK under the British Council Exchange Programme (HELS, 1989-90).  Currently, pursuing his dream of becoming a mathematician and exploring problems in Algebra, differential equations and Geometry.  

\vspace{24pt} 

\begin{center} 
\begin{tabular}{lcl}
\hline 
Dr Ramachandra Bhat,  & \hspace{50pt} & \textit{Alumni of}\\
115, 12$^{th}$ B Main Road, &  & 1. R V High School, Itgi, \\ 
6th Block, Rajajinagar, &  &   Taluk: Siddapur (Uttara Kannada), \\
Bengaluru - 560 010. &  & Karnataka, INDIA.  \\
Karnataka, INDIA. &  & 2. Sarada Vilas College, Mysore.\\
Mobile: +91-9902461175. &  & 3. Yuvaraja's College, Mysore. \\ 
Email: $rs\_bhat@yahoo.com $ &  & 4. Indian Institute of Technology, Bombay.\\
\hline 
\end{tabular} 
\end{center}   
 
\vspace{24pt}  
\begin{center}
\fbox{\textbf{Mathematicians are BORN to SOLVE Problems}} 
\end{center} 

\newpage 

\section*{Sample Data Tables}  

These tables are provided for the reference of article reviewer.  This could be deleted prior to publication, after the article review.  

\begin{center} 
\begin{longtable}{||c|c|c|c|c||}
  \multicolumn{5}{c} 
     {\textbf{Table 1.2: Collatz Sequence Data: n=(4P+3), P=(8x+1), x-Even}}\\ [5pt] 
     \hline 
  \multicolumn{5}{||c||}{PSO Computation: Collatz sequence terms \{$Col^a(n)$\} for n = (32x+7). } \\
  \hline 
  \multicolumn{1}{||c|}{$n_a$} 
  &\multicolumn{1}{|c|}{Formula} 
  &\multicolumn{1}{|c|}{Term value}
  &\multicolumn{1}{|c|}{Type}
  &\multicolumn{1}{|c||}{Condition}   \\  
  \hline 
\endfirsthead  
\multicolumn{5}{c}{\textbf{Table 1.2: The example data table ...} (\textit{Continued})}\\[5pt] 
\hline 
 \multicolumn{1}{||c|}{$n_a$} 
  &\multicolumn{1}{|c|}{Formula} 
  &\multicolumn{1}{|c|}{Term value}
  &\multicolumn{1}{|c|}{Type}
  &\multicolumn{1}{|c||}{Condition}   \\  
  \hline 
\endhead 
\hline 
\multicolumn{5}{r}{\small\itshape{(continued on the next page)}}\\
\endfoot 
\hline 
\multicolumn{5}{r}{\small \itshape{(This is just an example.  Complete compilation to be available as depository.}} 
\endlastfoot 
  \multicolumn{5}{||c||}{       } \\ 
\multicolumn{5}{||c||}{\textbf{Table 1.2.1: P=(8x+1); n=(32x+7); x - Even}} \\ 
\hline 
$n_1$ & \{3*(32x+7)+1\}       & (96x + 22) & Even  & - \\ 
$n_2$ & (96x + 22)$\div$ 2    & (48x + 11) & Odd   & - \\ 
$n_3$ & \{3*(48x+11)+1\}      & (144x + 34)& Even  & - \\ 
$n_4$ & (144x + 34) $\div$ 2  & (72x + 17) & Odd   & - \\ 
$n_5$ & \{3*(72x + 17) + 1\}  & (216x + 52)& Even  & - \\ 
$n_6$ & (216x + 52) $\div$ 2  & (108x + 26)& Even  & - \\ 
$n_7$ & (108x + 26) $\div$ 2  & (54x + 13) & Odd   & - \\ 
$n_8$ & \{3*(54x + 13)+1\}    & (162x + 40)& Even  & - \\ 
$n_9$ & (162x + 40) $\div$ 2  & (81x + 20) & Even  & x is Even\\
$n_{10}$ & (81x + 20) $\div$ 2  & (40x + 10) + ($\frac{x}{2}$) & Even & ($\frac{x}{2}$) - Even \\
$n_{11}$ & (81x + 20) $\div$ 4  & (20x + 5) + ($\frac{x}{4}$) & \{$<$ n\} & ($\frac{x}{2}$) - Even \\
\hline
\multicolumn{5}{||c||}{                } \\
\multicolumn{5}{||c||}{Result: Since $n_{11} < n $, \textbf{PSO = 11}, if \{x - Even and ($\frac{x}{2}$) - Even\}} \\ 
\multicolumn{5}{||c||}{\textbf{x=(8q+0) or (8q+4), for all q = \{0, 1, 2, 3, ...(All +ve integers)\}}}\\ 
\hline 
\multicolumn{5}{||c||}{                 } \\ 
\multicolumn{5}{||c||}{\textbf{Table 1.2.2: P = (8x+1); n = (32x+7); x - Even}} \\  
\hline 
$n_{10}$ & (81x + 20) $\div$ 2 & (40x + 10) + ($\frac{x}{2}$) & Odd  &  ($\frac{x}{2}$) - Odd \\
\hline 
$n_{11}$ & \{3*$\frac{81x + 20}{2}$ +1\} & ($\frac{243x+62}{2}$) &  Even & Even only \\ 
$n_{12}$ & (243x + 62) $\div$ 4  & (60x + 15) + ($\frac{3x+2}{4}$) & Even & ($\frac{3x+2}{4}$) - Odd \\ 
$n_{13}$ & (243x + 62) $\div$ 8  & (30x + 7) + ($\frac{3x+6}{8}$) & \{$<$ n\} & ($\frac{3x+2}{4}$) - Odd \\
\hline 
\multicolumn{5}{||c||}{                } \\
\multicolumn{5}{||c||}{Result: Since $n_{13} < n$, \textbf{PSO = 13}, if \{x-Even and ($\frac{3x+2}{4}$) - Odd\}} \\ 
\multicolumn{5}{||c||}{\textbf{x = (8q+6), for all q = \{0, 1, 2, 3, ...(All positive integers) \}}}\\  
\hline 
\multicolumn{5}{||c||}{                    } \\
\multicolumn{5}{||c||}{\textbf{Table 1.2.3: P = (8x+1); n = (32x+7); x - Even}} \\  
\hline 
$n_{12}$ & (243x + 62) $\div$ 4  & (60x + 15) + ($\frac{3x+2}{4}$) & Odd &  ($\frac{3x+2}{4}$) - Even \\ 
\hline 
$n_{13}$ & \{3*($\frac{243x + 62}{4}$) + 1\} & ($\frac{729x + 190}{4}$) & Even & Even only \\
$n_{14}$ & (729x + 190)$\div$8 & (91x + 23) + ($\frac{x+6}{8}$) & Even & ($\frac{x+6}{8}$)-Odd \\
$n_{15}$ & (729x + 190)$\div$16 & (45x + 11) + ($\frac{9x+14}{16}$) & Even & ($\frac{9x+14}{16}$)-Odd \\
$n_{16}$ & (729x + 190)$\div$32 & (22x + 5) + ($\frac{25x+30}{32}$) & \{$<$ n\} & ($\frac{9x+14}{16}$)-Odd \\
\hline 
\multicolumn{5}{||c||}{        } \\  
\multicolumn{5}{||c||}{Result: Since $n_{16} < n$, \textbf{PSO = 16}, if \{x-Even and ($\frac{9x+14}{16}$)-Odd\}} \\ 
\multicolumn{5}{||c||}{\textbf{x = (8q+2), q=(8r+2) \& (8r+6), for all r = \{0, 1, 2, 3, ...\}}}\\  
\hline 
\multicolumn{5}{||c||}{\textbf{Now onwards, it is all x=(8q+2), in these Tables}}  \\
\hline 
\multicolumn{5}{||c||}{            } \\
\multicolumn{5}{||c||}{\textbf{Table 1.2.4: P = (8x+1); n = (32x+7); x - Even}} \\  
\hline 
$n_{14}$ & (729x + 190)$\div$8 & (91x + 23) + ($\frac{x+6}{8}$) & Odd & ($\frac{x+6}{8}$)-Even \\
\hline 
$n_{15}$ & \{3*($\frac{729x + 190}{8}$) + 1\} & ($\frac{2187x + 578}{8}$) & Even & Even only \\
$n_{16}$ & (2187x + 578)$\div$16 & (136x + 36) + ($\frac{11x+2}{16}$) & Even & ($\frac{11x+2}{16}$)-Even \\
$n_{17}$ & (2187x + 578)$\div$32 & (68x + 18) + ($\frac{11x+2}{32}$) & Even & ($\frac{11x+2}{32}$)-Even \\
$n_{18}$ & (2187x + 578)$\div$64 & (34x + 9) + ($\frac{11x+2}{64}$) & Even & ($\frac{11x+2}{64}$)-Odd \\
$n_{19}$ & (2187x + 578)$\div$128 & (17x + 4) + ($\frac{11x+66}{128}$) & \{$<$ n\} & ($\frac{11x+2}{64}$)-Odd \\
\hline 
\multicolumn{5}{||c||}{        } \\  
\multicolumn{5}{||c||}{Result: Since $n_{19} < n$, \textbf{PSO = 19}, if \{x-Even and ($\frac{11x+2}{64}$)-Odd\}} \\ 
\multicolumn{5}{||c||}{\textbf{x = (8q+2), q=(8r+7), for all r=\{1, 3, 5, ...(All odd integers)\}}}\\  
\hline 
\multicolumn{5}{||c||}{                } \\
\multicolumn{5}{||c||}{\textbf{Table 1.2.5: P = (8x+1); n = (32x+7); x - Even}} \\  
\hline 
$n_{15}$ & (729x + 190)$\div$16 & (45x + 11) + ($\frac{9x+14}{16}$) & Odd & ($\frac{9x+14}{16}$)-Even \\
\hline 
$n_{16}$ & \{3*($\frac{729x + 190}{16}$) + 1\} & ($\frac{2187x + 586}{16}$) & Even & Even only \\
$n_{17}$ & (2187x + 586)$\div$32 & (68x + 18) + ($\frac{11x+10}{32}$) & Even & ($\frac{11x+10}{32}$)-Even \\
$n_{18}$ & (2187x + 586)$\div$64 & (34x + 9) + ($\frac{11x+10}{64}$) & Even & ($\frac{11x+10}{64}$)-Odd \\
$n_{19}$ & (2187x + 586)$\div$128 & (17x + 4) + ($\frac{11x+74}{128}$) & \{$<$ n\} & ($\frac{11x+10}{64}$)-Odd \\
\hline 
\multicolumn{5}{||c||}{        } \\  
\multicolumn{5}{||c||}{Result: Since $n_{19} < n$, \textbf{PSO = 19}, if \{x-Even and ($\frac{11x+10}{64}$)-Odd\}} \\ 
\multicolumn{5}{||c||}{\textbf{x = (8q+2), q=(8r+4), for all r=\{1, 3, 5, ...(All odd integers\}}}\\  
\hline 
\multicolumn{5}{||c||}{                } \\
\multicolumn{5}{||c||}{\textbf{Table 1.2.6: P = (8x+1); n = (32x+7); x - Even}} \\  
\hline 
$n_{16}$ & (2187x + 578)$\div$16 & (136x + 36) + ($\frac{11x+2}{16}$) & Odd & ($\frac{11x+2}{16}$)-Odd \\
\hline 
$n_{17}$ & \{3*($\frac{2187x + 578}{16}$) + 1\} & ($\frac{6561x + 1750}{16}$) & Even & Even only \\
$n_{18}$ & (6561x + 1750)$\div$32 & (205x + 54) + ($\frac{x+22}{32}$) & Even & ($\frac{x+22}{32}$)-Even \\
$n_{19}$ & (6561x + 1750)$\div$64 & (102x + 27) + ($\frac{33x+22}{64}$) & Even & ($\frac{33x+22}{64}$)-Odd \\
$n_{20}$ & (6561x + 1750)$\div$128 & (51x + 13) + ($\frac{33x+86}{128}$) & Even & ($\frac{33x+86}{128}$)-Odd \\
$n_{21}$ & (6561x + 1750)$\div$256 & (25x + 6) + ($\frac{161x+214}{256}$) & \{$<$ n\} & ($\frac{33x+86}{128}$)-Odd \\
\hline 
\multicolumn{5}{||c||}{        } \\  
\multicolumn{5}{||c||}{Result: Since $n_{21} < n$, \textbf{PSO = 21}, if \{x-Even and ($\frac{33x+86}{128}$)-Odd\}} \\ 
\multicolumn{5}{||c||}{\textbf{x = (8q+2), q=(8r+5), r=(8s+3, 8s+7), for all s=\{0, 1, 2, 3, ...\}}}\\  
\hline 
\multicolumn{5}{||c||}{                } \\
\multicolumn{5}{||c||}{\textbf{Table 1.2.7: P = (8x+1); n = (32x+7); x - Even}} \\  
\hline 
$n_{17}$ & (2187x + 578)$\div$32 & (68x + 18) + ($\frac{11x+2}{32}$) & Odd & ($\frac{11x+2}{32}$)-Odd \\
\hline 
$n_{18}$ & \{3*($\frac{2187 + 578}{32}$) + 1\} & ($\frac{6561x + 1766}{32}$) & Even & Even only \\
$n_{19}$ & (6561x + 1766)$\div$64 & (102x + 27) + ($\frac{33x+38}{64}$) & Even & ($\frac{33x+38}{64}$)-Odd \\
$n_{20}$ & (6561x + 1766)$\div$128 & (51x + 13) + ($\frac{33x+102}{128}$) & Even & ($\frac{33x+102}{128}$)-Odd \\
$n_{21}$ & (6561x + 1766)$\div$256 & (25x + 6) + ($\frac{161x+230}{256}$) & \{$<$ n\} & ($\frac{33x+102}{128}$)-Odd \\
\hline 
\multicolumn{5}{||c||}{        } \\  
\multicolumn{5}{||c||}{Result: Since $n_{21} < n$, \textbf{PSO = 21}, if \{x-Even and ($\frac{33x+102}{128}$)-Odd\}} \\ 
\multicolumn{5}{||c||}{\textbf{x = (8q+2), q=(8r+3), r=(8s+3, 8s+7), for all s=\{0, 1, 2, 3, ...\}}}\\  
\hline 
\multicolumn{5}{||c||}{                } \\
\multicolumn{5}{||c||}{\textbf{Table 1.2.8: P = (8x+1); n = (32x+7); x - Even}} \\  
\hline 
$n_{18}$ & (2187x + 578)$\div$64 & (34x + 9) + ($\frac{11x+2}{64}$) & Odd & ($\frac{11x+2}{64}$)-Even \\
\hline 
$n_{19}$ & \{3*($\frac{2187x + 578}{64}$) + 1\} & ($\frac{6561x + 1798}{64}$) & Even & Even only \\ 
$n_{20}$ & (6561x + 1798)$\div$128 & (51x + 14) + ($\frac{33x+6}{128}$) & Even & ($\frac{33x+6}{128}$)-Even \\
$n_{21}$ & (6561x + 1798)$\div$256 & (25x + 7) + ($\frac{161x+6}{256}$) & \{$<$ n\} & ($\frac{33x+6}{128}$)-Even \\
\hline 
\multicolumn{5}{||c||}{        } \\  
\multicolumn{5}{||c||}{Result: Since $n_{21} < n$, \textbf{PSO = 21}, if \{x-Even and ($\frac{33x+6}{128}$)-Even\}} \\ 
\multicolumn{5}{||c||}{\textbf{x = (8q+2), q=(8r+7), r=(8s+2, 8s+6), for all s=\{0, 1, 2, 3, ...\}}}\\  
\hline 
\multicolumn{5}{||c||}{                } \\
\multicolumn{5}{||c||}{\textbf{Table 1.2.9: P = (8x+1); n = (32x+7); x - Even}} \\  
\hline 
$n_{17}$ & (2187x + 586)$\div$32 & (68x + 18) + ($\frac{11x+10}{32}$) & Odd & ($\frac{11x+10}{32}$)-Odd \\
\hline 
$n_{18}$ & \{3*($\frac{2187x + 586}{32}$) + 1\} & ($\frac{6561x + 1790}{32}$) & Even & Even only \\ 
$n_{19}$ & (6561x + 1790)$\div$64 & (102x + 27) + ($\frac{33x+62}{64}$) & Even & ($\frac{33x+62}{64}$)-Odd \\ 
$n_{20}$ & (6561x + 1790)$\div$128 & (51x + 13) + ($\frac{33x+126}{128}$) & Even & ($\frac{33x+126}{128}$)-Odd \\ 
$n_{21}$ & (6561x + 1790)$\div$256 & (25x + 6) + ($\frac{161x+254}{256}$) & \{$<$ n\} & ($\frac{33x+126}{128}$)-Odd \\ 
\hline  
\multicolumn{5}{||c||}{                     } \\  
\multicolumn{5}{||c||}{Result: Since $n_{21} < n$, \textbf{PSO = 21}, if \{x-Even and ($\frac{33x+126}{128}$)-Odd\}} \\ 
\multicolumn{5}{||c||}{\textbf{x = (8q+2), q=(8r+0), r=(8s+3, 8s+7), for all s=\{0, 1, 2, 3, ...\}}}\\  
\hline  
\multicolumn{5}{||c||}{                } \\ 
\multicolumn{5}{||c||}{\textbf{Table 1.2.10: P = (8x+1); n = (32x+7); x - Even}} \\  
\hline 
$n_{18}$ & (2187x + 586)$\div$64 & (34x + 9) + ($\frac{11x+10}{64}$) & Odd & ($\frac{11x+10}{64}$)-Even \\
\hline 
$n_{19}$ & \{3*($\frac{2187x + 586}{64}$) + 1\} & ($\frac{6561x + 1822}{64}$) & Even & Even only \\ 
$n_{20}$ & (6561x + 1822)$\div$128 & (51x + 14) + ($\frac{33x+30}{128}$) & Even & ($\frac{33x+30}{128}$)-Even \\ 
$n_{21}$ & (6561x + 1822)$\div$256 & (25x + 7) + ($\frac{161x+30}{256}$) & \{$<$ n\} & ($\frac{33x+30}{128}$)-Even \\ 
\hline 
\multicolumn{5}{||c||}{        } \\  
\multicolumn{5}{||c||}{Result: Since $n_{21} < n$, \textbf{PSO = 21}, if \{x-Even and ($\frac{33x+30}{128}$)-Even\}} \\ 
\multicolumn{5}{||c||}{\textbf{x = (8q+2), q=(8r+4), r=(8s+2, 8s+6), for all s=\{0, 1, 2, 3, ...\}}}\\  
\hline 
\multicolumn{5}{||c||}{                } \\
\multicolumn{5}{||c||}{\textbf{Table 1.2.11: P = (8x+1); n = (32x+7); x - Even}} \\  
\hline 
$n_{18}$ & (6561x + 1750)$\div$32 & (205x + 54) + ($\frac{x+22}{32}$) & Odd & ($\frac{x+22}{32}$)-Odd \\
\hline 
$n_{19}$ & \{3*($\frac{6561x + 1750}{32}$) + 1\} & ($\frac{19683x + 5282}{32}$) & Even & Even only \\ 
$n_{20}$ & (19683x + 5282)$\div$64 & (307x + 82) + ($\frac{35x+34}{64}$) & Even & ($\frac{35x+34}{64}$)-Even \\ 
$n_{21}$ & (19683x + 5282)$\div$128 & (153x + 41) + ($\frac{99x+34}{128}$) & Even & ($\frac{99x+34}{128}$)-Odd \\ 
$n_{22}$ & (19683x + 5282)$\div$256 & (76x + 20) + ($\frac{227x+162}{256}$) & Even & ($\frac{227x+162}{256}$)-Even \\ 
$n_{23}$ & (19683x + 5282)$\div$512 & (38x + 10) + ($\frac{227x+162}{512}$) & Even & ($\frac{227x+162}{512}$)-Even \\ 
$n_{24}$ & (19683x + 5282)$\div$1024 & (19x + 5) + ($\frac{227x+162}{1024}$) & \{$<$ n\} & ($\frac{227x+162}{512}$)-Even \\ 
\hline 
\multicolumn{5}{||c||}{        } \\  
\multicolumn{5}{||c||}{Result: Since $n_{24} < n$, \textbf{PSO = 24}, if \{x-Even and ($\frac{227x+162}{512}$)-Even\}} \\ 
\multicolumn{5}{||c||}{\textbf{x = (8q+2), q=(8r+1), r=(8s+6), s=\{1, 3, 5, ...(All odd integers)\}}}\\  
\hline 
\multicolumn{5}{||c||}{                } \\
\multicolumn{5}{||c||}{\textbf{Table 1.2.12: P = (8x+1); n = (32x+7); x - Even}} \\ 
\hline 
$n_{19}$ & (6561x + 1750)$\div$64 & (102x + 27) + ($\frac{33x+22}{64}$) & Odd & ($\frac{33x+22}{64}$)-Even \\ 
\hline 
$n_{20}$ & \{3*($\frac{6561x + 1750}{64}$) + 1\} & ($\frac{19683x + 5314}{64}$) & Even & Even only \\ 
$n_{21}$ & (19683x + 5314)$\div$128 & (153x + 41) + ($\frac{99x+66}{128}$) & Even & ($\frac{99x+66}{128}$)-Odd \\ 
$n_{22}$ & (19683x + 5314)$\div$256 & (76x + 20) + ($\frac{227x+194}{256}$) & Even & ($\frac{227x+194}{256}$)-Even \\ 
$n_{23}$ & (19683x + 5314)$\div$512 & (38x + 10) + ($\frac{227x+194}{512}$) & Even & ($\frac{227x+194}{512}$)-Even \\ 
$n_{24}$ & (19683x + 5314)$\div$1024 & (19x + 5) + ($\frac{227x+194}{1024}$) & \{$<$ n\} & ($\frac{227x+194}{512}$)-Even \\ 
\hline 
\multicolumn{5}{||c||}{        } \\  
\multicolumn{5}{||c||}{Result: Since $n_{24} < n$, \textbf{PSO = 24}, if \{x-Even and ($\frac{227+194}{512}$)-Even\}} \\ 
\multicolumn{5}{||c||}{\textbf{x = (8q+2), q=(8r+5), r=(8s+0), for all s=\{1, 3, 5, ...(All odd integers)\}}}\\  
\hline 
\multicolumn{5}{||c||}{                } \\
\multicolumn{5}{||c||}{\textbf{Table 1.2.13: P = (8x+1); n = (32x+7); x - Even}} \\  
\hline 
$n_{20}$ & (6561x + 1750)$\div$128 & (51x + 13) + ($\frac{33x+86}{128}$) & Odd & ($\frac{33x+86}{128}$)-Even \\ 
\hline 
$n_{21}$ & \{3*($\frac{6561x + 1750}{128}$) + 1\} & ($\frac{19683x + 5378}{128}$) & Even & Even only \\ 
$n_{22}$ & (19683x + 5378)$\div$256 & (76x + 21) + ($\frac{227x+2}{256}$) & Even & ($\frac{227x+2}{256}$)-Odd \\ 
$n_{23}$ & (19683x + 5378)$\div$512 & (38x + 10) + ($\frac{227x+258}{512}$) & Even & ($\frac{227x+258}{512}$)-Even \\ 
$n_{24}$ & (19683x + 5378)$\div$1024 & (19x + 5) + ($\frac{227x+258}{1024}$) & \{$<$ n\} & ($\frac{227x+258}{512}$)-Even \\ 
\hline 
\multicolumn{5}{||c||}{        } \\  
\multicolumn{5}{||c||}{Result: Since $n_{24} < n$, \textbf{PSO = 24}, if \{x-Even and ($\frac{227x+258}{128}$)-Odd\}} \\ 
\multicolumn{5}{||c||}{\textbf{x = (8q+2), q=(8r+5), r=(8s+5), for all s=\{1, 3, 5, ...(All odd integers)\}}}\\  
\hline 
\multicolumn{5}{||c||}{                           } \\
\multicolumn{5}{||c||}{\textbf{Table 1.2.14: P = (8x+1); n = (32x+7); x - Even}} \\  
\hline  
$n_{19}$ & (6561x + 1766)$\div$64 & (102x + 27) + ($\frac{33x+38}{64}$) & Odd & ($\frac{33x+38}{64}$)-Even \\  
\hline  
$n_{20}$ & \{3*($\frac{6561x + 1766}{64}$) + 1\} & ($\frac{19683x + 5362}{64}$) & Even & Even only \\ 
$n_{21}$ & (19683x + 5362)$\div$128 & (153x + 41) + ($\frac{99x+114}{128}$) & Even & ($\frac{99x+114}{128}$)-Odd \\ 
$n_{22}$ & (19683x + 5362)$\div$256 & (76x + 20) + ($\frac{227x+242}{256}$) & Even & ($\frac{227x+242}{256}$)-Even \\  
$n_{23}$ & (19683x + 5362)$\div$512 & (38x + 10) + ($\frac{227x+242}{512}$) & Even & ($\frac{227x+242}{512}$)-Even \\ 
$n_{24}$ & (19683x + 5362)$\div$1024 & (19x + 5) + ($\frac{227x+242}{1024}$) & \{$<$ n\} & ($\frac{227x+242}{512}$)-Even \\ 
\hline 
\multicolumn{5}{||c||}{        } \\  
\multicolumn{5}{||c||}{Result: Since $n_{24} < n$, \textbf{PSO = 24}, if \{x-Even and ($\frac{227x+242}{512}$)-Even\}} \\ 
\multicolumn{5}{||c||}{\textbf{x = (8q+2), q=(8r+3), r=(8s+0), for all s=\{0, 2, 4, ...(All even integers)\}}}\\  
\hline 
\multicolumn{5}{||c||}{                } \\ 
\multicolumn{5}{||c||}{\textbf{Table 1.2.15: P = (8x+1); n = (32x+7); x - Even}} \\  
\hline 
$n_{20}$ & (6561x + 1766)$\div$128 & (51x + 13) + ($\frac{33x+102}{128}$) & Odd & ($\frac{33x+102}{128}$)-Even \\
\hline  
$n_{21}$ & \{3*($\frac{6561x + 1766}{128}$) + 1\} & ($\frac{19683x + 5426}{128}$) & Even & Even only \\  
$n_{22}$ & (19683x + 5426)$\div$256 & (76x + 21) + ($\frac{227x+50}{256}$) & Even & ($\frac{227x+50}{256}$)-Odd \\  
$n_{23}$ & (19683x + 5426)$\div$512 & (38x + 10) + ($\frac{227x+306}{512}$) & Even & ($\frac{227x+306}{512}$)-Even \\ 
$n_{24}$ & (19683x + 5426)$\div$1024 & (19x + 5) + ($\frac{227x+306}{1024}$) & \{$<$ n\} & ($\frac{227x+306}{512}$)-Even \\ 
\hline  
\multicolumn{5}{||c||}{               } \\  
\multicolumn{5}{||c||}{Result: Since $n_{24} < n$, \textbf{PSO = 24}, if \{x-Even and ($\frac{227x+306}{512}$)-Odd\}} \\ 
\multicolumn{5}{||c||}{\textbf{x = (8q+2), q=(8r+3), r=(8s+5), for all s=\{0, 2, 4, ...(All even integers)\}}}\\  
\hline 
\multicolumn{5}{||c||}{                } \\
\multicolumn{5}{||c||}{\textbf{Table 1.2.16: P=(8x+1); n=(32x+7); x-Even}} \\  
\hline 
$n_{20}$ & (6561x + 1798)$\div$128 & (51x + 14) + ($\frac{33x+6}{128}$) & Odd & ($\frac{33x+6}{128}$)-Odd \\
\hline  
$n_{21}$ & \{3*($\frac{6561x + 1798}{128}$) + 1\} & ($\frac{19683x + 5522}{128}$) & Even & Even only \\  
$n_{22}$ & (19683x + 5522)$\div$256 & (76x + 21) + ($\frac{227x+146}{256}$) & Even & ($\frac{227x+146}{256}$)-Odd \\  
$n_{23}$ & (19683x + 5522)$\div$512 & (38x + 10) + ($\frac{227x+402}{512}$) & Even & ($\frac{227x+402}{512}$)-Even \\ 
$n_{24}$ & (19683x + 5522)$\div$1024 & (19x + 5) + ($\frac{227x+402}{1024}$) & \{$<$ n\} & ($\frac{227x+402}{512}$)-Even \\ 
\hline  
\multicolumn{5}{||c||}{                      } \\  
\multicolumn{5}{||c||}{Result: Since $n_{24} < n$, \textbf{PSO = 24}, if \{x-Even and ($\frac{227x+402}{512}$)-Even\}} \\ 
\multicolumn{5}{||c||}{\textbf{i.e.,q=(8r+7), r=(8s+4), for all s = \{0, 2, 4,... (All even integers)\}}}\\  
\hline 
\multicolumn{5}{||c||}{                } \\
\multicolumn{5}{||c||}{\textbf{Table 1.2.17: P = (8x+1); n = (32x+7); x - Even}} \\  
\hline 
$n_{19}$ & (6561x + 1790)$\div$64 & (102x + 27) + ($\frac{33x+62}{64}$) & Odd & ($\frac{33x+62}{64}$)-Even \\ 
\hline 
$n_{20}$ & \{3*($\frac{6561x + 1790}{64}$)+ 1\} & ($\frac{19683x + 5434}{64}$)& Even & Even only \\
$n_{21}$ & (19683x + 5434)$\div$128  & (153x + 42) + ($\frac{99x+58}{128}$) & Even & ($\frac{99x+58}{128}$)-Even \\ 
$n_{22}$ & (19683x + 5434)$\div$256  & (76x + 21) + ($\frac{227x+58}{256}$) & Even & ($\frac{227x+58}{256}$)-Odd \\ 
$n_{23}$ & (19683x + 5434)$\div$512 & (38x + 10) + ($\frac{227x+314}{512}$) & Even & ($\frac{227x+314}{512}$)-Even \\ 
$n_{24}$ & (19683x + 5434)$\div$1024 &(19x + 5) + ($\frac{227x+314}{1024}$) & \{$<$ n\} & ($\frac{227x+314}{512}$)-Even \\ 
\hline  
\multicolumn{5}{||c||}{                     } \\  
\multicolumn{5}{||c||}{Result: Since $n_{24} < n$, \textbf{PSO = 24}, if \{x-Even and ($\frac{227x+314}{512}$)-Even\}} \\ 
\multicolumn{5}{||c||}{\textbf{x = (8q+2), q=(8r+x), r=(8s+x), for all s=\{0, 1, 2, 3, ...\}}}\\  
\hline  
\multicolumn{5}{||c||}{                } \\
\multicolumn{5}{||c||}{\textbf{Table 1.2.18: P = (8x+1); n = (32x+7); x - Even}} \\  
\hline 
$n_{20}$ & (6561x + 1790)$\div$128 & (51x + 13) + ($\frac{33x+126}{128}$) & Odd & ($\frac{33x+126}{128}$)-Even \\ 
\hline 
$n_{21}$ & \{3*($\frac{6561x + 1790}{128}$)+ 1\} & ($\frac{19683x + 5498}{128}$)& Even & Even only \\
$n_{22}$ & (19683x + 5498)$\div$256  & (76x + 21) + ($\frac{227x+122}{256}$) & Even & ($\frac{227x+122}{256}$)-Odd \\ 
$n_{23}$ & (19683x + 5498)$\div$512 & (38x + 10) + ($\frac{227x+378}{512}$) & Even & ($\frac{227x+378}{512}$)-Even \\ 
$n_{24}$ & (19683x + 5498)$\div$1024 &(19x + 5) + ($\frac{227x+378}{1024}$) & \{$<$ n\} & ($\frac{227x+378}{512}$)-Even \\ 
\hline  
\multicolumn{5}{||c||}{                     } \\  
\multicolumn{5}{||c||}{Result: Since $n_{24} < n$, \textbf{PSO = 24}, if \{x-Even and ($\frac{227x+378}{512}$)-Even\}} \\ 
\multicolumn{5}{||c||}{\textbf{x = (8q+2), q=(8r+0), r=(8s+1), for all s=\{0, 2, 4, ...(All Even only)\}}}\\  
\hline  
\multicolumn{5}{||c||}{                   } \\ 
\multicolumn{5}{||c||}{\textbf{Table 1.2.19: P = (8x+1); n = (32x+7); x - Even}} \\  
\hline 
$n_{20}$ & (6561x + 1822)$\div$128 & (51x + 14) + ($\frac{33x+30}{128}$) & Odd & ($\frac{33x+30}{128}$)-Odd \\ 
\hline 
$n_{21}$ & \{3*($\frac{6561x + 1822}{128}$)+ 1\} & ($\frac{19683x + 5594}{128}$)& Even & Even only \\
$n_{22}$ & (19683x + 5594)$\div$256  & (76x + 21) + ($\frac{227x+218}{256}$) & Even & ($\frac{227x+218}{256}$)-Odd \\ 
$n_{23}$ & (19683x + 5594)$\div$512 & (38x + 10) + ($\frac{227x+474}{512}$) & Even & ($\frac{227x+474}{512}$)-Even \\ 
$n_{24}$ & (19683x + 5594)$\div$1024 &(19x + 5) + ($\frac{227x+378}{1024}$) & \{$<$ n\} & ($\frac{227x+474}{512}$)-Even \\ 
\hline   
\multicolumn{5}{||c||}{        } \\  
\multicolumn{5}{||c||}{Result: Since $n_{24} < n$, \textbf{PSO = 24}, if \{x-Even and ($\frac{227x+474}{512}$)-Even\}} \\ 
\multicolumn{5}{||c||}{\textbf{x = (8q+2), q=(8r+4), r=(8s+0), for all s=\{0, 2, 4, ... (All even only)\}}}\\  
\hline 
\multicolumn{5}{||c||}{                } \\
\multicolumn{5}{||c||}{\textbf{Table 1.2.20: P = (8x+1); n = (32x+7); x - Even}} \\  
\hline 
$n_{20}$ & (19683x + 5282)$\div$64 & (307x + 82) + ($\frac{35x+34}{64}$) & Odd & ($\frac{35x+34}{64}$)-Odd \\ 
\hline 
$n_{21}$ & \{3*($\frac{19683x + 5282}{64}$)+ 1\} & ($\frac{59049x + 15910}{64}$)& Even & Even only \\ 
$n_{22}$ & (59049x + 15910)$\div$128  & (461x + 124) + ($\frac{41x+38}{128}$) & Even & ($\frac{41x+38}{128}$)-Even \\ 
$n_{23}$ & (59049x + 15910)$\div$256  & (230x + 62) + ($\frac{169x+38}{256}$) & Even & ($\frac{169x+38}{256}$)-Even \\ 
$n_{24}$ & (59049x + 15910)$\div$512 & (115x + 31) + ($\frac{169x+38}{512}$) & Even & ($\frac{169x+38}{512}$)-Odd \\ 
$n_{25}$ & (59049x + 15910)$\div$1024 &(57x + 15) + ($\frac{681x+550}{1024}$) & Even & ($\frac{681x+550}{1024}$)-Odd \\ 
$n_{26}$ & (59049x + 15910)$\div$2048 &(28x + 7) + ($\frac{1705x+1574}{2048}$) & \{$<$ n\} & ($\frac{681x+550}{1024}$)-Odd \\ 
\hline  
\multicolumn{5}{||c||}{        } \\  
\multicolumn{5}{||c||}{Result: Since $n_{26} < n$, \textbf{PSO = 26}, if \{x-Even and ($\frac{681x+550}{1024}$)-odd\}, i.e.,} \\ 
\multicolumn{5}{||c||}{\textbf{x = (8q+2), q=(8r+1), r=(8s+5), s=(8t+2, 8t+6), for all t=\{0, 1, 2, 3, ...\}}}\\  
\hline 
\end{longtable} 
\end{center} 
 

\begin{thebibliography}{1} 

\bibitem{Tao} %1
\textit{The Notorious Collatz Conjecture}, Presentation note by Prof. Terence Tao, UCLA, 2019. https://terrytao.files.wordpress.com/2020/02/collatz.pdf   

\bibitem{Peter} %2 
Are we near a solution to the 3x+1 problem?  
www.occampress.com/solution.pdf

\end{thebibliography}
\end{document}